%% file: GarnierO_-_Uniqueness_of_roots_up_to_conjugacy_in_circular_and_hosohedral-type_Garside_groups.tex
\documentclass[11pt,twoside]{amsart}
\usepackage{amsfonts}
\usepackage{hyperref}
\usepackage{etoolbox}
\usepackage{amsmath}
\usepackage{amssymb}
\usepackage{amsthm}
\usepackage{dsfont}
\usepackage{pifont}
\usepackage{mathrsfs}
\usepackage{graphicx}
\usepackage{enumerate}
\usepackage[all]{xy}
\usepackage{geometry}
\usepackage{amsfonts}
\usepackage{fancyhdr}
\usepackage{tikz}
\usepackage{tikz-cd}
\usetikzlibrary{arrows}
\usepackage[section]{algorithm}
\usepackage{algorithmic}

\usepackage{amsaddr}

\usepackage[backend=biber,style=alphabetic,sorting=nty,maxnames=5]{biblatex}
\include{pream.tex}%toutes les commandes y sont enregistrees

\theoremstyle{plain}
\newtheorem{prop}{Proposition}[section]
\newtheorem{prop-def}[prop]{Proposition-Definition}
\newtheorem{lem}[prop]{Lemma}
\newtheorem{theo}[prop]{Theorem}
\newtheorem{cor}[prop]{Corollary}

\newtheorem*{theo*}{Theorem}
\newtheorem*{cor*}{Corollary}
\newtheorem*{prop*}{Proposition}
\newtheorem*{fact*}{Fact}

\theoremstyle{definition}
\newtheorem{definition}[prop]{Definition}

\theoremstyle{remark}
\newtheorem{rem}[prop]{Remark}
\newtheorem{exemple}[prop]{Example}

 \normalbaselines
\title[Uniqueness of roots in circular and hosohedral-type Garside groups]{Uniqueness of roots up to conjugacy in circular and hosohedral-type Garside groups}
\author{Owen Garnier}
\address{LAMFA, Université de Picardie Jules Verne, CNRS UMR 7352,\\ 33, rue Saint-Leu, 80000, Amiens, France.}
\email{o.garnier@u-picardie.fr}
\date{\today}

\begin{document}

\begin{abstract}
We consider a particular class of Garside groups, which we call circular groups. We mainly prove that roots are unique up to conjugacy in circular groups. This allows us to completely classify these groups up to isomorphism. As a consequence, we obtain the uniqueness of roots up to conjugacy in complex braid groups of rank 2.

We also consider a generalization of circular groups, called hosohedral-type groups. These groups are defined using circular groups, and a procedure called the $\Delta$-product, which we study in generality. We also study the uniqueness of roots up to conjugacy in hosohedral-type groups.\end{abstract}

\newcommand{\pr}{\mathrm{pr}}
\renewcommand{\sl}{\mathrm{sl}}

\maketitle
\tableofcontents
\addtocontents{toc}{\protect\setcounter{tocdepth}{0}}
\section*{Introduction}
The aim of this article is to study two particular classes of Garside monoids, called circular and hosohedral-type monoids, along with their associated Garside groups. Circular monoids were first defined in the original article of Dehornoy and Paris where the authors introduce the notion of Garside group. Hosohedral-type monoids first appeared in the work of Matthieu Picantin (\cite{thespicantin}) and later on in that of Mireille Soergel (\cite{soergelsystolic}) because of the particular properties of their lattices of simples.

Let $m,\ell$ be two positive integers. The circular monoid $M(m,\ell)$ is defined by a monoid presentation, given by a set of $m$ ordered generators, endowed with relations stating that any product of $\ell$ consecutive generators should be equal (Definition \ref{def:circular_monoid}). The circular group $G(m,\ell)$ is then defined as the enveloping group of $M(m,\ell)$. Circular monoids are Garside monoids, with Garside element the product of $\ell$ consecutive generators.

Our main motivation is to study roots in circular groups using their Garside structure, in particular, the associated solution to the conjugacy problem. We reduce the situation to the study of two classes of elements, namely \emph{rigid} and \emph{periodic} elements (Proposition \ref{prop:rigid_or_periodic}).

Generally speaking, let $(M,\Delta)$ be a Garside monoid, and let $p,q$ be integers. An element $\rho$ in the associated Garside group $G(M)$ is said to be $(p,q)$-periodic if $\rho^p=\Delta^q$. We give a complete description of periodic elements in circular groups, up to conjugacy.

\begin{prop*}(Proposition \ref{prop:irreducible_periodic_circulaire} and Corollary \ref{cor:center_circular_group}) Let $m,\ell$ be two positive integers. An element $x\in G(m,\ell)$ is periodic if and only if it is conjugate to a power of a product of either $m$ or $\ell$ consecutive generators. In particular, for a given $p$ and $q$, all $(p,q)$-periodic elements are conjugate. Moreover, an element $x\in G(m,\ell)$ is periodic if and only if it has a central power.
\end{prop*}

As a consequence, an isomorphism between two circular groups must preserve periodic elements. From this we deduce the classification of circular groups up to isomorphism.

\begin{cor*}(Corollary \ref{cor:classification_circular_groups}) Let $m,\ell,m',\ell'$ be four positive integers. The groups $G(m,\ell)$ and $G(m',\ell')$ are isomorphic if and only if one of the following holds
\begin{enumerate}[-]
\item $1\in \{m,\ell\}$ and $1\in \{m',\ell'\}$. In this case, $G(m,\ell)\simeq G(m',\ell')\simeq \Z$.
\item $(m',\ell')\in \{(m,\ell),(\ell,m)\}$.
\end{enumerate}
\end{cor*}

Since periodic elements are (by definition) roots of powers of $\Delta$, the above proposition is a first result of uniqueness of roots up to conjugacy in circular groups. By also studying the conjugacy of rigid elements, we obtain the following stronger result.

\begin{theo*}(Theorem \ref{theo:roots_conjugate_circular}) Let $m, \ell$ be positive integers. If $\alpha,\beta\in G(m,\ell)$ are such that $\alpha^n=\beta^n$ for some nonzero integer $n$, then $\alpha$ and $\beta$ are conjugate in $G(m,\ell)$.
\end{theo*}

We then apply this result to complex braid groups of rank $2$. Recall that, if $W\subset \GL_n(\C)$ is a finite group generated by (pseudo-)reflections, then one can consider the braid group $B(W)$, defined as the fundamental group of the regular orbit space associated to $W$ (see \cite[Section 2.B]{bmr}). These groups generalize Artin groups of spherical type.

The question of uniqueness of roots up to conjugacy in braid groups was first studied by Juan Gonz\'alez-Meneses in \cite{uniquenessrootbraids1} , where he proved that roots were unique up to conjugacy in the Artin groups of type $A$ (i.e. the usual braid groups). His results were later expanded in \cite{uniquenessrootbraids2} to the Artin groups of type $B$. To our knowledge, no other results are known regarding this question. The uniqueness of roots up to conjugacy is conjectured to hold at least for every spherical Artin group (cf. \cite[Conjecture X.3.10]{ddgkm}). 

The above theorem on uniqueness of roots up to conjugacy in circular groups, along with the fact that every complex braid group of rank $2$ is isomorphic to a circular group (see \cite[Theorem 1 and Theorem 2]{bannai}), gives the uniqueness of roots up to conjugacy in complex braid groups of rank $2$.

\begin{theo*}(Theorem \ref{theo:roots_conjugate_rank2}) Let $W$ be an irreducible complex reflection group of rank $2$, and let $B(W)$ be its braid group. If $\alpha,\beta\in B(W)$ are such that $\alpha^n=\beta^n$ for some nonzero integer $n$, then $\alpha$ and $\beta$ are conjugate in $B(W)$.
\end{theo*}

This theorem gives preliminary evidences that the uniqueness of roots up to conjugacy may actually hold in all complex braid groups.

In the last section, we study a generalization of circular monoids, called hosohedral-type monoids. These monoid appear as a particular case of a general construction involving Garside monoids, called the \emph{$\Delta$-product}.

Let $(M_1,\Delta_1),\ldots,(M_h,\Delta_h)$ be a family of Garside monoids. We define the $\Delta$-product of $(M_1,\Delta_1)*_\Delta\cdots *_\Delta(M_h,\Delta_h)$ as the quotient of the free product of the $M_i$ by the relations stating that all the $\Delta_i$ are equal to each other. This monoid is again a Garside monoid (with Garside element the image of any $\Delta_i$), in which each one of the $M_i$ embeds.

\begin{prop*}(Proposition \ref{prop:left-weighted_facto_in_delta_product}) Let $(M_1,\Delta_1),\ldots,(M_h,\Delta_h)$ be homogeneous Garside monoids. Each Garside group $G(M_i)$ naturally embeds in the Garside group $G(M_1*_\Delta\cdots *_\Delta M_h)$. These embeddings respect left-weighted factorizations.
\end{prop*}

From this proposition, we deduce that the center of a $\Delta$-product of Garside monoids (with more than one factor) is cyclic and generated by a power of the Garside element (Proposition \ref{prop:centre_delta_product}). We also deduce that, for integers $p,q$, the $\Delta$-product $M_1*_\Delta\cdots *_\Delta M_h$ admits $(p,q)$-periodic elements if and only if one of the factors does (Proposition \ref{prop:periodic_delta_product}). 
%More precisely, the embeddings $G(M_i)\to G(M)$ preserve periodic elements, and every periodic element in $G(M)$ is conjugate to the image in $G(M)$ of some periodic element of some $G(M_i)$.

Note that two periodic elements coming from two different factors $G(M_i)$ and $G(M_j)$ may not be conjugate in the $\Delta$-product (Example \ref{ex:periodic_torus}).

We then apply these results to hosohedral-type monoids, which are defined as the $\Delta$-product of circular groups. In this case, there can be distinct conjugacy classes of $(p,q)$-periodic elements (see Example \ref{ex:periodic_torus}). We show that these are actually the only obstructions to the uniqueness of roots up to conjugacy.

\begin{theo*}(Theorem \ref{theo:roots_conjugate_hosohedral}) Let $M$ be a hosohedral-type monoid. If $\alpha,\beta\in G(M)$ are such that $\alpha^n=\beta^n$ for some nonzero integer $n$, then we either have
\begin{enumerate}[$\bullet$]
\item $\alpha$ and $\beta$ are conjugate.
\item $\alpha$ and $\beta$ are nonconjugate periodic elements of $G(M)$.
\end{enumerate}
\end{theo*}

\subsection*{Acknowledgments.} This work is part of my PhD thesis, done under the supervision of Pr. Ivan Marin. I thank him for his precious advice, especially in Sections \ref{sec:homology_circular} and \ref{sec:isomorphisms}. I would also like to thank Mireille Soergel and Igor Haladjian for stimulating discussions.
\addtocontents{toc}{\protect\setcounter{tocdepth}{2}}

\section{Preliminaries on Garside monoids}
Throughout this article, the gcd (resp. lcm) of two integers $p$ and $q$ will be denoted by $p\wedge q$ (resp. $p \vee q$).

Our main arguments rely on the study of super-summit sets and periodic elements. For this we need some reminders on Garside monoids. Our main reference is \cite[Section I.2 and Chapter VIII]{ddgkm}.

\subsection{Definitions, normal form}

We start by considering a monoid $M$. Throughout this paper, we will always assume that $M$ is \nit{homogeneous}. That is, there is some monoid morphism $\ell:M\to (\Z_{\geqslant 0},+)$ such that elements of positive length generate $M$. This condition is far from minimal when considering Garside monoids, but it is sufficient here as both circular monoids and hosohedral-type monoids are homogeneous. The homogeneity condition implies in particular that the set of invertible elements of $M$ is trivial. We also assume that $M$ is \nit{cancellative}. That is, every equality of the form $abc=ab'c$ in $M$ implies $b=b'$. Under these assumptions, we define two partial orders $\preceq$ and $\succeq$ on $M$ by
\[a\preceq b\Leftrightarrow \exists c\in M~|~ ac=b \text{~and~} b\succeq a\Leftrightarrow \exists c\in M~|~ b=ca.\]
These partial order are called \nit{left-divisibility} and \nit{right-divisibility}, respectively. A nontrivial element $a\in M$ is an \nit{atom} if it admits no nontrivial (right- or left-)divisors other than itself. We say that an element $\Delta\in M$ is \nit{balanced} if its sets of left- and right-divisors are equal, we then simply call this set the divisors of $\Delta$.

\begin{definition}\cite[Definition I.2.1]{ddgkm} \newline Let $M$ be a homogeneous cancellative monoid, such that the posets $(M,\preceq)$ and $(M,\succeq)$ are both lattices. A \nit{Garside element} in $M$ is a balanced element $\Delta$ whose set of divisors $S$ is finite and generates $M$. We call $S$ the set of \nit{simple elements} (associated to $\Delta$) and we say that $(M,\Delta)$ is a \nit{homogeneous Garside monoid}.
\end{definition}

The assumption that $(M,\preceq)$ and $(M,\succeq)$ are both lattices means that (left- and right-) lcms and gcds of two elements always exist. We denote by $a\wedge b$ (resp. $a\vee b$) the left-gcd (resp. right-lcm) of two elements $a,b$. If need be, we will denote by $\wedge_R$ the right-gcd and by $\vee_L$ the left-lcm.

Let $(M,\Delta)$ be a homogeneous Garside monoid, and let $s\in S$ be a simple element. By cancellativity, there is a unique simple $\bbar{s}$ such that $s\bbar{s}=\Delta$. We call $\bbar{s}$ the \nit{left-complement} of $s$ in $\Delta$. Likewise, the \nit{right-complement} $s^*$ of $s$ in $\Delta$ is defined by $s^*s=\Delta$. 

Defining $\phi(s):=\bbar{\bbar{s}}$ for any simple element $s$ yields an automorphism of $M$, which permutes the simple elements (in particular, it has finite order). We call $\phi$ the \nit{Garside automorphism} of $(M,\Delta)$. By definition, we have $s\Delta=\Delta\phi(s)$ for any simple element $s$.

One of the main features of a Garside monoid is that it gives rise to a convenient solution to the word problem, given by the notion of greedy word.
\begin{definition}\cite[Corollary V.1.54]{ddgkm}\newline Let $(M,\Delta)$ be a homogeneous Garside monoid. A word $st$ of length 2 in $S$ is called \nit{greedy} if $s=(st)\wedge \Delta$, or equivalently if $\bbar{s}$ and $t$ are left-coprime. In general, a word $s_1\cdots s_r$ in $S$ is called \nit{greedy} if each subword $s_is_{i+1}$ is greedy for $i\in \intv{1,r-1}$.
\end{definition}

\begin{lem}\label{gr2}
Let $s,t$ be two simple elements in a homogeneous Garside monoid $(M,\Delta)$, and let $u:=t\wedge \bbar{s}$. The greedy normal form of the product $st$ is given by $s't'$ where
\[s'=su\in S \text{ and }t'=u^{-1}t\in S.\]
\end{lem}

\begin{proof}
This is an application of the invariance of left-gcd under multiplication, as we have $su=s(t\wedge \bbar{s})=(st)\wedge (s\bbar{s})=x\wedge \Delta$.
\end{proof}

An immediate induction on this lemma shows that any element in a homogeneous Garside monoid $M$ admits a unique decomposition as a greedy word in the simple elements. The lemma also gives a practical way to compute the greedy decomposition of any element $x$ of $M$, by the means of successively computing greedy decompositions of pairs of elements of $S$ (see for instance \cite[Algorithm III.1.52]{ddgkm}).

The definition of a homogeneous Garside monoid implies the Ore condition, thus a homogeneous Garside monoid $M$  always embeds in its enveloping group $G(M)$, which can be conveniently described as a group of fractions. The greedy decomposition in $M$ carries on into a complete description of the elements of $G(M)$.

\spa\begin{prop}\cite[Proposition I.2.4]{ddgkm}\newline
Let $(M,\Delta)$ be a homogeneous Garside monoid. Every element $x\in G(M)$ admits a unique decomposition of the form $x=\Delta^kb$ with $k\in \Z,b\in M$ and $\Delta\not\preceq b$.
\end{prop}
In particular, this proposition gives a solution to the word problem in $G(M)$, given by computing the above decomposition, and then the greedy decomposition of $b$ in $M$. This decomposition of an element $x\in G(M)$ is the \nit{left-weighted factorization} of $x$.

\subsection{Conjugacy and periodic elements}
In this section, we fix $(M,\Delta)$ a homogeneous Garside monoid, and $G(M)$ its enveloping group. The Garside monoid $M$ provides a solution to the conjugacy problem in $G(M)$, given by the computation of super-summit sets. These super-summit sets also allow for the computation of centralizers in $G(M)$. 

Let $x\in G(M)$ with left-weighted factorization $x=\Delta^ks_1\cdots s_r$. The \nit{infimum} (resp. the \nit{supremum}) of $x$ is defined by $\inf(x):=k$ (resp. $\sup(x)=k+r$).

\begin{definition}\cite[Defintion VIII.2.12]{ddgkm}\newline
Let $x\in G(M)$. The conjugacy class of $x$ in $G(M)$ admits a well-defined subset $\SSS(x)$, on which each one of $\inf$ and $\sup$ takes a constant value. Furthermore, for every conjugate $x'$ of $x$ in $G(M)$, we have
\[\inf(x)\leqslant \inf(\SSS(x))\text{~and~}\sup(x')\geqslant \sup(\SSS(x)).\]
The set $\SSS(x)$ is called the \nit{super-summit set} of $x$. 
\end{definition}

Since the set $S$ is finite, there is a finite number of elements of $G(M)$ with given $\inf$ and $\sup$. In particular, the super-summit set of any element of $G(M)$ is always finite. As the super-summit set of an element depends only on its conjugacy class, computing this set gives a solution to the conjugacy problem. To compute super-summit sets, one uses the cycling and decycling operations.

Let $x\in G(M)$ with left-weighted factorization $x=\Delta^ks_1\cdots s_r$. The \nit{initial factor} (resp. \nit{final factor}) of $x$ is defined as $\init(x):=\phi^{-k}(s_1)$ (resp. $\fin(x):=s_r$). The \nit{cycling} (resp. \nit{decycling}) of $x$ is then defined as
\[\cyc(x)=x^{\init(x)}=x^{\phi^{-k}(s_1)}~(\text{resp. }\dec(x)=x^{\fin(x)^{-1}}=x^{s_r^{-1}}).\]

\begin{prop}\cite[Proposition VII.2.16]{ddgkm}\label{prop:reach_sss}\newline Let $x\in G(M)$. One can go from $x$ to an element of $\SSS(x)$ by a finite sequence of cycling, followed by a finite sequence of decycling.
\end{prop}

This proposition can be used to compute an element of $\SSS(x)$ starting from $x\in G(M)$. The whole super-summit set of $x\in G(M)$ can then be encoded in a so-called conjugacy graph.

\begin{definition}\cite[Section 3]{centralizergarside}\newline 
Let $x\in G(M)$, the \nit{conjugacy graph} $\CG(x)$ of $x$ is an oriented graph, defined as follows.
\begin{enumerate}[$\bullet$]
\item The object set of $\CG(x)$ is $\SSS(x)$.
\item An arrow $a\to b$ in $\CG(x)$ is given by a simple element $s\in M$ such that $a^s=b$ and that no left-divisor $t$ of $s$ is such that $a^t\in \CG(x)$.
\end{enumerate}
\end{definition}

By \cite[Lemma VIII.2.19]{ddgkm}, the conjugacy graph of any element $x$ is connected. In other words, two elements of $\SSS(x)$ are always conjugate by a sequence of simple elements (and their inverses) such that the result at each step remains in $\SSS(x)$.

Let $x\in G(M)$. By definition of the conjugacy graph, a non-oriented path from $x$ to itself in $\CG(x)$ induces a word in the simple elements (and their inverses) representing an element of the centralizer $C_{G(M)}(x)$. By \cite[Theorem 3.4]{centralizergarside}, the centralizer $C_{G(M)}(x)$ of $x$ is generated by the images in $G(M)$ of any generating set of paths from $x$ to itself in $\CG(x)$.

In our study of circular and hosohedral-type monoids, we are going to consider two classes of elements showing two extreme behavior with respect to the cycling and decycling operations.

\begin{definition}\cite[Definition 3.1]{rigidity}\newline
An element $x\in G(M)$, with left-weighted factorization $x=\Delta^ks_1\cdots s_r$, is said to be \nit{rigid} if $r=0$ or if the word $\fin(x)\init(x)$ is greedy. 
\end{definition}

If $x=\Delta^ks_1\cdots s_r\in G(M)$ is rigid (with $r>0$), the left-weighted factorizations of $\cyc(x)$ and $\dec(x)$ are given by
\[\cyc(x)=\Delta^k s_2\cdots s_r\phi^{-k}(s_1)\text{~and~}\dec(x)=\Delta^k \phi^{k}(s_r)s_1\cdots s_{r-1}.\]
In particular, $\inf(\cyc(x))=\inf(x)=\inf(\dec(x))$ and $\sup(\cyc(x))=\sup(x)=\sup(\dec(x))$. As both $\cyc(x)$ and $\dec(x)$ are also rigid, we obtain that a rigid element always belongs to its super-summit set by Proposition \ref{prop:reach_sss}.

Rigid elements also exhibit a strict behavior with respect to powers. This behavior is useful in relating the super-summit set of a rigid element to that of its powers.

\begin{lem}\label{lem:cycling_power_rigid}
Let $(M,\Delta)$ be a Garside monoid, and let $\alpha,\beta\in G(M)$ be rigid elements. For any positive integer $n$, we have $\cyc(\alpha)^n=\cyc(\alpha^n)$, $\dec(\alpha)^n=\dec(\alpha^n)$ and $\phi(\alpha^n)=\phi(\alpha)^n$. Moreover, if $\alpha^n=\beta^n$ for some positive integer $n$, then $\alpha=\beta$.
\end{lem}
\begin{proof}
The result is trivial if $n=1$, we assume that $n\geqslant 2$ from now on. Let $\alpha=\Delta^ks_1\cdots s_r$ be the left-weighted factorization of $\alpha$. Since $\alpha$ is rigid, the left-weighted factorization of $x:=\alpha^n$ is 
\[\alpha^n=\Delta^{nk}\phi^{(n-1)k}(s_1\cdots s_r)\cdots \phi^k(s_1\cdots s_r)s_1\cdots s_r.\]
We note that the initial factor (resp. the final factor) of $\alpha^n$ is the same as that of $\alpha$, hence the result on $\cyc(\alpha^n)$ and $\dec(\alpha^n)$. The result on $\phi(\alpha^n)$ is an obvious consequence of $\phi(x)=x^\Delta$ for $x\in G(M)$.

Moreover, we have $\inf(\alpha^n)=n\inf(\alpha)$ and $\sup(\alpha^n)=n\sup(\alpha)$. The left-weighted factorization of $\alpha$ can be recovered using that of $\alpha^n$ by taking $\Delta^{\frac{\inf(\alpha^n)}{n}}$, followed by the last $\frac{\sup(\alpha^n)}{n}$ terms of the left-weighted factorization of $\alpha^n$. Since this depends only on the left-weighted factorization of $\alpha^n$, we obtain that $\alpha^n=\beta^n$ implies $\alpha=\beta$ if $\beta$ is another rigid element.
\end{proof}

The other class of elements we will be interested in is that of periodic elements.

\begin{definition}\cite[Definition V.3.2]{ddgkm}\newline Let $(M,\Delta)$ be a homogeneous Garside monoid, and let $p,q$ be two integers. An element $\rho\in G(M)$ is said to be $(p,q)$\nit{-periodic} if $\rho^p=\Delta^q$.
\end{definition}

By definition, any root and any power of a periodic element in $G(M)$ is again a periodic element. Moreover, it is obvious by definition that a $(p,q)$-periodic element is also $(np,nq)$-periodic for any integer $n\geqslant 1$. One can show that the converse is also true:

\begin{prop}\cite[Proposition VIII.3.31 and Proposition VIII.3.34]{ddgkm}\label{prop:periodic_elements}\newline 
Let $p,q$ be two integers. Let $d:=p\wedge q$, with $p=dp'$ and $q=dq'$. If $p'=1$, then any $(p,q)$-periodic element in $G(M)$ is conjugate to $\Delta^{q'}$. If $p'\neq 1$, then any $(p,q)$-periodic element in $G(M)$ is conjugate to a $(p',q')$-periodic element of the form $\Delta^k s$, where $s\in S$.
\end{prop}

In particular, this proposition shows that the super-summit set of a periodic element contains only elements of the form $\Delta^k s$ or $\Delta^k$ (as such elements clearly lie in their own super-summit set).

The definition of rigid element is a priori not group-theoretic and depends on the Garside monoid $M$. On the contrary, periodic elements can sometimes be characterized by a solely group-theoretic property. If $(M,\Delta)$ is a homogeneous Garside monoid such that the center of $G(M)$ is cyclic and generated by a power of $\Delta$, then an element $\rho\in G(M)$ is periodic if and only if it admits a central power. In this case, we can try and compare two groups by comparing their respective periodic elements.

\begin{definition}
Let $(M,\Delta)$ be a homogeneous Garside monoid. We say that a periodic element $\rho\in G(M)$ is \nit{irreducible} if it admits no roots in $G(M)$ other than itself.
\end{definition}

Since $(M,\Delta)$ is assumed to be homogeneous here, we obtain that any periodic element in $G(M)$ is always a power of some irreducible periodic elements. Let $(M,\Delta),(M',\Delta')$ be two homogeneous Garside monoids, such that the center of $G(M)$ (resp. of $G(M')$) is cyclic and generated by some power of $\Delta$ (resp. $\Delta'$). An isomorphism $G(M)\to G(M')$ sends a generator of $Z(G(M))$ to a generator of $Z(G(M'))$. Therefore, it must map (irreducible) periodic elements of $G(M)$ to (irreducible) periodic elements of $G(M')$.

\section{Circular monoids}\label{sec:circular}

In this section we define circular monoids by their (monoid) presentations. These monoids already appeared in \cite[Example 5]{dehpar}, where the authors showed that they were Garside monoids. Here we propose an in-depth study of their Garside properties. 

\subsection{Definition, first properties}

Let $m,\ell$ be two positive integers that we fix throughout this section. Let also $\{a_0,\ldots,a_{m-1}\}$ be an alphabet. For $i\in \Z$ and $p\in \Z_{\geqslant 1}$, we define $s(i,p)$ as the word
\[s(i,p):=\prod_{k=i'}^p a_k=a_{i'}a_{i'+1}\cdots a_{i'+p-1}\]
where $i'$ is the remainder in the Euclidean division of $i$ by $m$, and with the convention that, for $j\geqslant 1$, $a_j:=a_{j'}$ where $j'$ is the remainder in the Euclidean division of $j$ by $m$. We also define $s(i,0)$ to be the empty word for all $i\in \Z$.

\begin{definition}\label{def:circular_monoid}
Let $m,\ell$ be two positive integers. The \nit{circular monoid} $M(m,\ell)$ is defined by the monoid presentation 
\[M(m,\ell):=\left\langle a_0,\ldots,a_{m-1}~|~ \forall i\intv{0,m-1}, s(i,\ell)=s(i+1,\ell)\right\rangle^+.\]
The enveloping group $G(m,\ell)$ of $M(m,\ell)$ is called a \nit{circular group}.
\end{definition}

From now on, we assimilate the word $s(i,p)$ ($i\in \Z$, $p\geqslant 0$) with its image in $M(m,\ell)$.

\begin{exemple}
The monoid $M(3,3)$ is given by $\langle a,b,c~|~abc=bca=cab\rangle^+$. The group $G(3,3)$ is the fundamental group of the complement of $3$ lines going through the origin in $\C^2$ (cf. \cite[Example 5]{dehpar}). The monoid $M(2,3)$ is given by $\langle s,t~|~sts=tst\rangle^+$. It is the Artin monoid of type $A_2$.
\end{exemple}

We note that the presentation of $M(m,\ell)$ is homogeneous. That is, the defining relations are equalities between words of the same length. The function sending an element of $M(m,\ell)$ to the length of any word representing it is then a length function, making $M(m,\ell)$ into a homogeneous monoid. As the only defining relations of $M(m,\ell)$ are between words of length $\ell$, two words of length less than $\ell$ cannot represent the same element of $M(m,\ell)$. In particular, for $0<p<\ell$ and $i,i'\in \Z$, we have $s(i,p)=s(i',p)$ in $M(m,\ell)$ if and only if $i\equiv i'[m]$. 

\begin{lem}\cite[Example 5]{dehpar}\newline
The monoid $M(m,\ell)$ is a homogeneous Garside monoid with Garside element $\Delta=s(0,\ell)$. Its simple elements are the $s(i,p)$ for $i\in \intv{0,m-1}$ and $p\in \intv{0,\ell}$. The Garside automorphism of $M(m,\ell)$ is given by $\phi(s(i,p))=s(i+\ell,p)$.
\end{lem}

The only statement which is not showed in \cite[Example 5]{dehpar} is the statement on the Garside automorphism, which comes from the fact that, for any simple element $s(i,p)$ of $M(m,\ell)$, we have
\[s(i,p)\Delta=s(i,p)s(i+p,\ell)=s(i,\ell+p)=s(i,\ell)s(i+\ell,p)=\Delta s(i+\ell,p).\]

Let $s(i,p)$ be a simple element of $M(m,\ell)$. We have $s(i,p)s(i+p,\ell-p)=s(i,\ell)=\Delta$, thus the right-complement (resp. left-complement) of $s(i,p)$ in $\Delta$ is given by
\[s(i,p)^*=s(i+p-\ell,\ell-p)~(\text{resp. }\bbar{s(i,p)}=s(i+p,\ell-p)).\]

The fact that any simple element different from $1$ and $\Delta$ admits a unique decomposition as a product of atoms has the following consequences.

\begin{lem}(Left-gcd of two simple elements)\label{lem:gcd_lcm_simple_circulaire}\newline
Let $s(i,p)$ and $s(i',p')$ be two simple elements of $M(m,\ell)$. 
\begin{enumerate}[(a)]
\item We have $s(i,p)\preceq s(i',p')$ if and only if $i=i'$ and $p\leqslant p'$, or if $p=0$, or if $p'=\ell$.
\item The left-gcd of $s(i,p)$ and $s(i',p')$ is given by
\[s(i,p)\wedge s(i',p')=\begin{cases} s(i,p)&\text{if }s(i,p)\preceq s(i',p') \\ s(i',p')&\text{if }s(i',p')\preceq s(i,p)\\ 1&\text{otherwise} \end{cases}\]
\item The right-lcm of $s(i,p)$ and $s(i',p')$ is given by
\[s(i,p)\vee s(i',p')=\begin{cases} s(i',p')&\text{if }s(i,p)\preceq s(i',p')\\ s(i,p)&\text{if }s(i',p')\preceq s(i,p)\\ \Delta&\text{otherwise} \end{cases}\]
\end{enumerate}
\end{lem}
\begin{proof}
$(a)$ The cases $p=0$ and $p'=\ell$ are immediate: they give $s(i,p)=1$ and $s(i',p')=\Delta$, respectively. If $i=i'$ and $p\leqslant p'$, we have $s(i,p)s(i+p,p'-p)=s(i',p')$ and $s(i,p)\preceq s(i',p')$. Conversely, suppose that $s(i,p)\preceq s(i',p')$. Since we assume that $p'<\ell$ and $0<p$, there is only one way to write both $s(i,p)$ and $s(i',p')$ as a product of atoms. Thus the assumption that $s(i,p)$ can be written as a prefix of some word expressing $s(i',p')$ implies that $s(i,p)$ is in fact the only prefix of length $p$ of $s(i',p')$. We obtain that $p\leqslant p'$ and $s(i',p)=s(i,p)$, which gives $i=i'$. \newline $(b)$ The first two cases are obvious. Suppose that we have $s(i,p)\not\preceq s(i',p')$ and $s(i',p')\not\preceq s(i,p)$. By the first point we have $p,p'\notin\{0,\ell\}$ and $i\neq i'$. Again by point $(a)$, a nontrivial common left-divisor $s(j,q)$ of $s(i,p)$ and $s(i',p')$ should be such that $j=i=i'$, which is impossible. We apply similar reasoning to prove point $(c)$.
\end{proof}

This very strict behavior of gcds and lcms allows for an easy description of greedy normal forms of a product of two simple elements in $M(m,\ell)$.

\begin{lem}(Greedy normal form of a product of two nontrivial simples)\label{lem:greedy_form_2_simples}\newline 
Let $s(i,p)$ and $s(i',p')$ be two simple elements of $M(m,\ell)$ with $0<p,p'<\ell$. The greedy normal form of $s(i,p)s(i',p')$ is given by
\[s(i,p)s(i',p')=\begin{cases} s(i,p)s(i',p')&\text{if }i+p\not\equiv i' [m], \\ 
s(i,p+p') &\text{if }i+p\equiv i'[m] \text{ and }p+p'<\ell,\\
\Delta &\text{if }i+p\equiv i'[m] \text{ and }p+p'=\ell,\\
\Delta s(i+\ell,p+p'-\ell)&\text{if }i+p\equiv i' [m] \text{ and }p+p'> \ell. \end{cases}\]
\end{lem}
\begin{proof}
We will apply Lemma \ref{gr2}. Since $0<p,p'<\ell$, Lemma \ref{lem:gcd_lcm_simple_circulaire} gives
\begin{align*}
\bbar{s(i,p)}\wedge s(i',p')&=s(i+p,\ell-p)\wedge s(i',p'),\\
&=\begin{cases} s(i+p,\ell-p)&\text{if }i+p\equiv i'[m] \text{ and }\ell-p\leqslant p,\\ s(i',p')&\text{if }i+p\equiv i'[m] \text{ and }p'\leqslant \ell-p\\ 1&\text{otherwise}.\end{cases}
\end{align*}
These three cases give the desired result.
\end{proof}

This particular description of greedy normal forms in circular monoids will induce a convenient description of super-summit sets in Section \ref{sec:conjugacy_circular}.

We finish this section by the study of the particular case where $m=\ell$. In this case, the group $G(m,m)$ is the fundamental group of the complement of $m$ lines through the origin in $\C^2$ (cf. \cite[Example 5]{dehpar}). Its presentation can be reinterpreted as a direct product.

\begin{lem}\label{lem:m=l}
Let $m$ be a positive integer, and let $F_{m-1}$ be a free group with $m-1$ generators. The group $G(m,m)$ is isomorphic to $\Z\times F_{m-1}$. This isomorphism identifies $\Delta$ with $(z,1)$, where $z$ is a generator of $\Z$.
\end{lem}
\begin{proof}
The result is obvious if $m\in \{1,2\}$. We denote by $x_0,\cdots,x_{m-2}$ the generators of $F_{m-1}$, and by $z$ a generator of $\Z$. We define $f:\Z\times F_{m-1}\to G(m,m)$ by
\[\begin{cases} f(z)=\Delta=s(0,m),\\ f(x_i)=a_i&\forall i\in \intv{0,m-2}. \end{cases}\]
This induces a well-defined morphism since $\Delta\in Z(G(m,m))$. Conversely, we define $g:G(m,m)\to \Z\times F_{m-1}$ by
\[\begin{cases} g(a_i)=x_i&\forall i\in \intv{1,m-2},\\ g(a_{m-1})=(x_0\cdots x_{m-2})^{-1}z.\end{cases}\]
This induces a well-defined morphism since, for all $i\in \intv{0,m-1}$, we have
\begin{align*}
g(s(i,m))&=g(s(i,m-1-i)a_{m-1}s(0,i))\\
&=x_i\cdots x_{m-2} (x_0\cdots x_{m-2})^{-1} z x_0\cdots x_{i-1}\\
&=x_i\cdots x_{m-2} (x_0\cdots x_{m-2})^{-1} x_0\cdots x_{i-1}z=z=g(s(0,m)).
\end{align*}
It is straightforward to check that $f$ and $g$ are inverses of one another.
\end{proof}

\subsection{Conjugacy in circular groups}\label{sec:conjugacy_circular}
Again we fix two positive integers $m,\ell$. In this section, we study separately the conjugacy of periodic and non-periodic elements in circular groups. Our first result shows that one can reduce to the case of periodic or rigid elements.

\begin{prop}\label{prop:rigid_or_periodic}
Let $x$ be an element of $G(m,\ell)$. If $x$ lies in its own super-summit set, then it is either rigid or periodic.
\end{prop}

\begin{proof}
First, if $\inf(x)=\sup(x)$, then we have $x=\Delta^k$ for some $k\in \Z$. In this case, $x$ is obviously both rigid and $(1,k)$-periodic.

Suppose now that $\sup(x)=\inf(x)+1$. We have $x=\Delta^k s(i,p)$ for some $i\in \intv{0,m-1}$ and $p\in \intv{1,\ell-1}$. The element $x$ is rigid if and only if the word $s(i,p)\phi^{-k}(s(i,p))=s(i,p)s(i-k\ell,p)$ is greedy. By Lemma \ref{lem:greedy_form_2_simples}, this is equivalent to $i+p\equiv i-k\ell[m]$, i.e. $k\ell+p\equiv 0[m]$. If $k\ell+p\equiv 0[m]$, then there is some integer $v$ with $-k\ell+mv=p$. In particular we have $\ell \wedge m|p$. We have
\begin{align*}
x^n&=\Delta^{kn}\phi^{k(n-1)}(s(i,p))\cdots \phi^k(s(i,p))s(i,p)\\
&=\Delta^{kn}s(i+k(n-1)\ell,p)\cdots s(i+k\ell,p)s(i,p)\\
&=\Delta^{kn}s(i+k(n-1)\ell,np)\\
&=\Delta^{kn+a}s(i+(k(n-1)+a)\ell,r),
\end{align*}
where $np=a\ell+r$ is the Euclidean division of $np$ by $\ell$. The remainder $r$ is $0$ if and only if $n$ is a multiple of $\frac{\ell\vee p}{p}=\frac{\ell}{\ell\wedge p}$. We obtain that $x$ is $(a,b)$-periodic, for
\[a=\frac{\ell}{\ell\wedge p}\text{~and~}b=k\frac{\ell\vee p}{p}+\frac{\ell\vee p}{\ell}=\frac{k\ell+p}{\ell\wedge p}=\frac{mv}{\ell\wedge p}.\]
Furthermore, $a$ and $b$ are coprime.

Lastly, suppose that $\sup(x)>\inf(x)+1$. The left-weighted factorization of $x$ is given by $\Delta^ks_1\cdots s_r$ with $r>1$. We claim that $x$ is rigid. Otherwise, the word $s_r\phi^{-k}(s_1)$ is not greedy. By Lemma \ref{lem:greedy_form_2_simples}, we either have that $s_r\phi^{-k}(s_1)$ is a simple element, or a product of the form $\Delta s$ where $s$ is a simple element different from $1$ and $\Delta$. In the first case, we have $\sup(\cyc(x))<\sup(x)$. In the second case, we have $\inf(\cyc(x))>\inf(x)$. In both cases, we have $x\notin\SSS(x)$.
\end{proof}

\subsubsection{Periodic elements}
The proof of Proposition \ref{prop:rigid_or_periodic} also gives the following result.
\begin{lem}\label{lem:elements_periodiques_circulaires}
Let $m,\ell$ be positive integers, and let $\Delta^ks(i,p)$ be in $G(m,\ell)$ with $0<p<\ell$. The element $\Delta^k s(i,p)$ is periodic if and only if $p+k\ell\equiv 0[m]$, in which case it is a $(\frac{\ell}{p\wedge \ell},\frac{mv}{p\wedge \ell})$-periodic element (where $p+k\ell=mv$).
\end{lem}

Now that we have a characterization of the elements of the super-summit sets of periodic elements, we can compute conjugacy graphs and centralizers. We distinguish two cases.

\begin{lem}\label{lem:sss_deltak}
Let $x=\Delta^k$ for some nonzero integer $k$. We have $\SSS(\Delta^k)=\{\Delta^k\}$. The centralizer of $\Delta^k$ in $G(m,\ell)$ is either $G(m,\ell)$ if $k\ell$ is a multiple of $m$, or cyclic and generated by $\Delta$ otherwise. 
\end{lem}
\begin{proof}
We have that $y\in G(m,\ell)$ lies in $\SSS(\Delta^k)$ only if $\inf(y)=k=\sup(y)$. The only element satisfying this is $\Delta^k$, which is conjugate to itself. Thus we have $\SSS(\Delta^k)=\{\Delta^k\}$. Now, let $s(j,q)$ be a simple element in $M(m,\ell)$. Since both $1$ and $\Delta$ conjugate $\Delta^k$ to itself, we can assume that $q\in \intv{1,m-1}$. We have 
\begin{align*}
s(j,q)^{-1} \Delta^k s(j,q)&=\bbar{s(j,q)} \Delta^{k-1}s(j,q)\\
&=s(j+q,\ell-q)\Delta^{k-1}s(j,q)\\
&=\Delta^{k-1} s(j+q+(k-1)\ell,\ell-q)s(j,q).
\end{align*}
In order for this element to lie in $\SSS(\Delta^k)$, the word $s(j+q+(k-1)\ell,\ell-q)s(j,q)$ must not be greedy. This is equivalent to $j+k\ell\equiv j[m]$. If $k\ell$ is a multiple of $m$, this is true for all $j\in \intv{0,m-1}$, and we obtain $s(j,q)^{-1}\Delta^ks(j,q)=\Delta^k$: the arrows from $\Delta^k$ to itself in $\CG(\Delta^k)$ are given by all the simple elements. Otherwise, $j+k\ell\equiv j[m]$ is never true for $j\in \intv{1,m-1}$ and the only arrows from $\Delta^k$ to itself in $\CG(\Delta^k)$ are given by $1$ and $\Delta$.
\end{proof}

\begin{lem}\label{lem:sss_periodic}
Let $x=\Delta^ks(i,p)$ be a periodic element in $M(m,\ell)$ with $p\in \intv{1,m-1}$. We have $\SSS(x)=\{\Delta^k s(n,p)~|~ n\in \intv{0,m-1}\}$. The centralizer of $\Delta^ks(0,p)$ in $G(m,\ell)$ is cyclic and generated by $s(p,m)$.
\end{lem}
\begin{proof}The assumption that $x$ is periodic is equivalent to $k\ell+p\equiv 0[m]$ by Lemma \ref{lem:elements_periodiques_circulaires}. Let $s(j,q)$ be a simple element in $M(m,\ell)$. We have 
\begin{align*}
x^{s(j,q)}&=s(j,q)^{-1}\Delta^k s(i,p) s(j,q)\\
&=\Delta^{k-1} s(j+q+(k-1)\ell,\ell-q)s(i,p)s(j,q).
\end{align*}
Again, in order for this to lie in $\SSS(x)$, we must have either $j+k\ell\equiv i[m]$ or $i+p\equiv j[m]$. Since $k\ell+p\equiv 0[m]$, those two assertions are equivalent. If they are satisfied, then we have
\[x^{s(j,q)}=\Delta^{k}s(j+q-p,p)=\Delta^ks(i+q,p).\]
In particular, $s(p,n)$ gives a conjugating element from $\Delta^ks(0,p)$ to $\Delta^ks(n,p)$ for $n\in \intv{0,m-1}$. Moreover, for $\Delta^k s(n,p)\in \SSS(x)$, the simples $s$ such that $(\Delta^k s(n,p))^s\in \SSS(x)$ are all divisible by $s(n+p,1)$. The conjugacy graph of $x$ is then given by
\[\xymatrixcolsep{3.5pc}\xymatrix{\Delta^ks(0,p)\ar[r]^-{s(p,1)}&\Delta^ks(1,p)\ar[r]^-{s(p+1,1)}&\cdots \ar[r]^-{s(p+m-2,1)}&\Delta^k s(m-1,p) \ar@/^2pc/[lll]^-{s(p+m-1,1)}},\]
and the centralizer of $\Delta^ks(0,p)$ is cyclic and generated by 
\[s(p,1)s(p+1,1)\cdots s(p+m-1,1)=s(p,m).\]
\end{proof}

We can use these two lemmas to determine the center of circular groups. Recall that the Garside automorphism $\phi$, corresponding to conjugacy by $\Delta$ on the right, sends a simple element $s(i,p)$ to $s(i+\ell,p)$. If $m\neq 1\neq \ell$, then the smallest trivial power of $\phi$ is $\phi^{\frac{m}{m\wedge \ell}}$, and $\Delta^{\frac{m}{m\wedge \ell}}$ is the smallest central power of $\Delta$ in $G(m,\ell)$.

\begin{cor}(Center of circular groups)\label{cor:center_circular_group}\newline Let $m,\ell$ be two positive integers. If $m=1$ or $\ell=1$, then $G(m,\ell)\simeq \Z$ is abelian. If $m=\ell=2$, then $G(m,\ell)=\Z^2$ is abelian. Otherwise $Z(G(m,\ell))$ is infinite cyclic and generated by $\Delta^{\frac{m}{m\wedge \ell}}$.
\end{cor}
\begin{proof}
If $m=1$, then $M(1,\ell)=\langle a_0 \rangle^+\simeq \Z_{\geqslant 0}$ (with Garside element $a_0^\ell$). If $\ell=1$, we have $M(m,1)=\langle a_0\rangle^+\simeq \Z_{\geqslant 0}$ (with Garside element $a_0$). If $m=\ell=2$, then $G(m,\ell)=\langle a_0,a_1~|~a_0a_1=a_1a_0\rangle=\Z^2$ (with Garside element $a_0a_1$).

If $m>1$, we distinguish several cases. First, we assume that $m$ does not divide $\ell$. By Lemma \ref{lem:sss_deltak}, the centralizer of $\Delta$ in $G(m,\ell)$ is cyclic and generated by $\Delta$. As the center $Z(G(m,\ell))$ is included in $C_{G(m,\ell)}(\Delta)$, we obtain that $Z(G(m,\ell))$ is cyclic and generated by the smallest central power of $\Delta$, which is $\Delta^{\frac{m}{m\wedge \ell}}$ since $\ell\neq 1$.

Now, if $m=\ell$, then Lemma \ref{lem:m=l} gives an isomorphism $G(m,\ell)\simeq \Z\times F_{m-1}$. The cases $m=\ell\in \{1,2\}$ have already been studied. If $m\geqslant 3$, then the center of $\Z\times F_{m-1}$ is $\Z\times \{1\}$, which is identified with $\langle \Delta \rangle=\langle \Delta^{\frac{m}{m\wedge \ell}}\rangle$.

Lastly, we assume that $mk=\ell$ for some integer $k>1$. The element $s(0,m)$ is $(k,1)$-periodic in $G(m,\ell)$ and, by Lemma \ref{lem:sss_periodic}, the centralizer of $s(0,m)$ in $G(m,\ell)$ is cyclic and generated by $s(0,m)$. Since $mk=\ell$, we have that $s(0,m)^{k}=\Delta=\Delta^{\frac{m}{m\wedge \ell}}$ is a central element. It remains to show that $s(0,m)$ admits no central power inferior to $k$. Let $1\leqslant r\leqslant k-1$. We have $s(0,m)^r=s(0,mr)$ and 
\begin{align*}
s(0,1)^{s(0,m)^r}&=s(0,1)^{s(0,mr)}\\
&=\Delta^{-1}s(mr-\ell,\ell-mr)s(0,1)s(0,mr)\\
&=\Delta^{-1}s(0,m(k-r))s(0,1)s(0,mr)
\end{align*}
As $m\neq 1$, this is a left-weighted factorization, in particular it is not equal to $s(0,1)$. We then obtain that the smallest central power of $s(0,m)$ is $s(0,m)^k=\Delta$, thus $Z(G(m,\ell))=\langle \Delta \rangle$ as claimed.
\end{proof}

\begin{prop}\label{prop:periodic_are_conjugate_circular}
Let $m,\ell$ be positive integers, and let $p,q$ be integers. Any two $(p,q)$-periodic elements in $G(m,\ell)$ are conjugate.
\end{prop}
\begin{proof}
First, by Proposition \ref{prop:periodic_elements}, we can assume that $p,q$ are coprime integers. If $p=1$, then $x$ and $y$ are both conjugate to $\Delta^q$. If $p>1$, we can assume (up to conjugacy) that $x,y$ are of the form $\Delta^ks(i,a)$ and $\Delta^{k'}s(i',a')$ with $0<a,a'<\ell$, respectively. By Lemma \ref{lem:elements_periodiques_circulaires}, there are two integers $v,v'$ with $mv=k\ell+a$ and $mv'=k'\ell+a'$. Since $k\ell+a$ (resp. $k'\ell+a'$) is the length of $x$ (resp. $y$) in $G(m,\ell)$, we have $v=v'$. By reducing modulo $\ell$, we obtain $a\equiv a'[\ell]$. Since we assume that $a,a'\in \intv{0,\ell-1}$, $a\equiv a'[\ell]$ implies $a=a'$. The equality $k\ell+a=k'\ell+a'$ then gives $k=k'$. By Lemma \ref{lem:sss_periodic}, we get that $x$ and $y$ are conjugate.
\end{proof}

\begin{prop}\label{prop:irreducible_periodic_circulaire}
Let $m,\ell$ be positive integers. Any periodic element in $G(m,\ell)$ is conjugate to a power of either $s(0,m)$ or $\Delta$. Moreover the irreducible periodic elements of $G(m,\ell)$ are given (up to conjugacy) by
\[\begin{cases} \{s(0,m)^{\pm 1}\}&\text{if }m|\ell,\\ \{\Delta^{\pm 1}\} &\text{if }\ell|m,\\ \{s(0,m)^{\pm 1},\Delta^{\pm 1}\}&\text{otherwise}.\end{cases}\]
\end{prop}
\begin{proof}
Let $\rho\in G(m,\ell)$ be a periodic element. By Lemma \ref{lem:sss_periodic}, we can assume up to conjugacy that $\rho=\Delta^k s(k\ell,p)$. If $p=0$, then $\rho=\Delta^k$ is a power of $\Delta$. If $p\neq 0$, then Lemma \ref{lem:elements_periodiques_circulaires} gives an integer $v$ with $mv=k\ell+p$, we then have $\rho=s(0,m)^v$. If $\rho$ is an irreducible periodic element of $G(m,\ell)$, then we have $\rho\in \{s(0,m)^{\pm 1},\Delta^{\pm 1}\}$ (up to conjugacy). It only remains to check whether or not $s(0,m)$ and $\Delta$ are indeed irreducible. If $m|\ell$ (resp. $\ell|m$), then we have $s(0,m)^{\frac{\ell}{m}}=\Delta$ (resp. $\Delta^{\frac{m}{\ell}}=s(0,m)$). Since there must be at least one conjugacy class of irreducible periodic elements in $G(m,\ell)$, we get the desired result if $m|\ell$ or $\ell|m$.

Assume now that neither $m|\ell$ nor $\ell|m$. A proper root of $s(0,m)$ in $G(m,\ell)$ must have the form $s(0,n)$ with $0<n<m$. By Lemma \ref{lem:elements_periodiques_circulaires}, such an element cannot be periodic, thus it cannot be a root of $s(0,m)$, which is then irreducible. The same reasoning applies to $\Delta$. 
\end{proof}

\subsubsection{Non-periodic elements}We now turn our attention to non-periodic elements. By Proposition \ref{prop:rigid_or_periodic}, such elements are exactly the conjugate of rigid elements in $G(m,\ell)$.

\begin{prop}\label{prop:conjugacy_graph_rigid_circular}
Let $x\in G(m,\ell)$ be a non-periodic element. The super-summit set of $x$ is made of rigid elements. Furthermore, the only arrows starting from an object $y$ of $\CG(x)$ are labeled by $\init(y)$ and $\bbar{\fin(y)}$.
\end{prop}
\begin{proof}
Let $y\in \SSS(x)$. Since $x$ is not periodic, $y$ is not periodic. It is then rigid by Proposition \ref{prop:rigid_or_periodic}. We then have $\sup(y)>\inf(y)$ and we can assume that the left-weighted factorization of $y$ is $\Delta^k s(i_1,p_1)\cdots s(i_r,p_r)$ with $r>0$. Since $y$ is rigid, we have $i_r+p_r\not\equiv i_1-k\ell [m]$ by Lemma \ref{lem:greedy_form_2_simples}. Let $s(j,q)$ be a simple element with $q\in \intv{1,m-1}$. We have
\begin{align*}
y^{s(j,q)}&=\Delta^{k-1}s(j+q+(k-1)\ell,\ell-q)s(i_1,p_1)\cdots s(i_r,p_r)s(j,q).
\end{align*}
In order for this to lie in $\SSS(x)$, we must have either $j+k\ell\equiv i_1[m]$ or $i_r+p_r\equiv j[m]$. Since $i_r+p_r+kl\not\equiv i_1[m]$, these cases are mutually exclusive.
\begin{enumerate}[-]
\item Assume that $j+k\ell\equiv i_1[m]$. By Lemma \ref{lem:greedy_form_2_simples}, the left-weighted factorization of $y^{s(j,q)}$ is given by
\[\begin{cases}\Delta^{k-1}s(j+q+(k-1)\ell,\ell-q+p_1)\cdots s(i_r,p_r)s(j,q)&\text{if }p_1<q,\\
\Delta^{k}s(i_2,p_2)\cdots s(i_r,p_r)s(j,q)&\text{if }p_1=q,\\
\Delta^{k}s(j+q+k\ell,p_1-q)s(i_2,p_2)\cdots s(i_r,p_r)s(j,q)&\text{if }p_1>q. \end{cases}\]
Thus, $y^{s(j,q)}\in \SSS(x)$ in this case if and only if $p_1=q$. We then have that $s(j,q)=s(i_1-k\ell,p)$ is the initial factor of $y$.
\item Assume that $i_r+p_r\equiv j[m]$. By Lemma \ref{lem:greedy_form_2_simples}, the left-weighted factorization of of $y^{s(j,q)}$ is given by
\[\begin{cases}\Delta^{k-1}s(j+q+k\ell-\ell,\ell-q)s(i_1,p_1)\cdots s(i_r,p_r+q)&\text{if }p_r+q<\ell,\\
\Delta^{k}\phi(s(j+q+k\ell-\ell,\ell-q)s(i_1,p_1)\cdots s(i_{r-1},p_{r-1}))&\text{if }p_r+q=\ell,\\
\Delta^{k}\phi(s(j+q+k\ell-\ell,\ell-q)s(i_1,p_1)\cdots s(i_{r-1},p_{r-1}))s(i_r+\ell,p_r+q-\ell)&\text{if }p_r+q>\ell. \end{cases}\]
Thus, $y^{s(j,q)}\in \SSS(x)$ in this case if and only if $p_r+q=\ell$. We then have that $s(j,q)=s(i_r+p_r,\ell-p_r)=\bbar{s(i_r,p_r)}$.
\end{enumerate}
\end{proof}

Let $x\in G(m,\ell)$ be a rigid element. The conjugation of $x$ by $\bbar{\fin(x)}$ is equal to $\phi(\dec(x))$. The last proposition then gives

\begin{cor}\label{cor:connecting_rigid_conjugates}
Let $x\in G(m,\ell)$ be a non-periodic element. One can go from any element of $\SSS(x)$ to any other by a finite sequence of cycling, decycling and application of the Garside automorphism.
\end{cor}

We can now state our main result on uniqueness of roots up to conjugacy.

\begin{theo}\label{theo:roots_conjugate_circular}(Uniqueness of roots up to conjugacy in circular groups)\newline
Let $m,\ell$ be two positive integers. If $\alpha,\beta\in G(m,\ell)$ are such that $\alpha ^n=\beta^n$ for some nonzero integer $n$, then $\alpha $ and $\beta$ are conjugate.
\end{theo}
\begin{proof}
First, if $\alpha $ is $(p,q)$-periodic for some integers $p$ and $q$. We have that $\alpha ^n$ is $(p,nq)$-periodic and that $\beta$ is also $(np,nq)$-periodic. The elements $\alpha $ and $\beta$ are then conjugate by Proposition \ref{prop:periodic_elements} and Proposition \ref{prop:periodic_are_conjugate_circular}.

Up to replacing $\alpha$ and $\beta$ with $\alpha^{-1}$ and $\beta^{-1}$, we can assume that $n>0$. Assume now that $\alpha $ is not periodic, we also have that $x:=\alpha ^n$ and $\beta$ are non periodic. Up to conjugacy, we can assume that $\alpha \in \SSS(\alpha )$. By Proposition \ref{prop:rigid_or_periodic}, we have that $\alpha $ is rigid. The element $x$ is then rigid as a power of the rigid element $\alpha $. Let now $c\in G(m,\ell)$ be so that $\beta^c\in \SSS(\beta)$. Since $\beta$ is not periodic, $\beta^c$ is rigid as well as $x^c=(\beta^c)^n$. We have $x,x^c\in \SSS(x)$. By Corollary \ref{cor:connecting_rigid_conjugates}, there is a finite sequence of cycling, decycling, and application of the Garside automorphism sending $x$ to $x^c$. By Lemma \ref{lem:cycling_power_rigid}, applying the same transformations to $\alpha $ gives a rigid element $\alpha '$ whose $n$-th power is $x^c$. Again by Lemma \ref{lem:cycling_power_rigid}, we have $\alpha '=\beta^c$ and thus $\alpha $ and $\beta$ are conjugate.
\end{proof}

\subsection{Some group theoretic properties}
\subsubsection{Homology of circular groups}\label{sec:homology_circular}
The homology of a Garside group can be studied using a particular complex introduced by Dehornoy and Lafont in \cite[Section 4]{dehlaf}. This complex is built using atoms and lcms in the underlying Garside monoid. The particular behavior of circular monoids with regards to lcms induces strong results on the associated complex.

Let us start by quickly recalling the definition of the Dehornoy-Lafont complex. We start by considering a homogeneous Garside monoid $(M,\Delta)$, with set of simples $S$ and set of atoms $A$. We fix an arbitrary strict linear ordering $<$ on $A$. For any $x\in M$, we define $\md(x)$ to be the $<$-minimal element of $A$ which right-divides $x$.

An $n$-cell is defined as an $n$-tuple $[\alpha_1,\ldots,\alpha_n]$ of atoms of $M$ such that $\alpha_1<\ldots<\alpha_n$, and $\alpha_i=\md(\alpha_i\vee_L\alpha_{i+1}\vee_L\cdots\vee_L \alpha_n)$ for all $i\in \intv{1,n}$. We denote by $\Xx_n$ the set of $n$-cells. The set $C_n$ of $n$-chains is then defined as the free $\Z G(M)$-module with basis the set of $n$-cells. Endowed with a convenient differential $\partial_n$ (which we will not define here), the complex $(C_n,\partial_n)_{n\in \N}$ is an exact resolution of the trivial $\Z G(M)$-module $\Z$. Since both $S$ and $A$ are finite, we have $\Xx_n=\varnothing$ for $n\geqslant |A|$, thus $(C_n,\partial_n)_{n\geqslant 0}$ is bounded above and below. 

\begin{lem}\label{lem:cellules_dehlaf_circular}
Let $m,\ell$ be two positive integers, and let $M:=M(m,\ell)$. We have $\Xx_0=\{[\varnothing]\}$, $\Xx_1=\{[a_0],\ldots,[a_{m-1}]\}$, $\Xx_2=\{[a_0,a_i]~|~i\in \intv{1,m-1}\}$ and $\Xx_n=\varnothing$ for $n\geqslant 3$.
\end{lem}
\begin{proof}
The statements on $\Xx_0$ and $\Xx_1$ are straightforward. By definition, an $n$-tuple $[\alpha_1,\ldots,\alpha_n]$ is an $n$-cell if and only if $[\alpha_2,\ldots,\alpha_n]$ is an $n-1$-cell and $\alpha_1=\md(\alpha_2\vee_L\cdots\vee_L\alpha_n)$.

If $n=2$, then we get that a couple $[a_i,a_j]$ is a $2$-cell if and only if $a_i=\md(a_i\vee_L a_j)$. Since $a_i\neq a_j$ by assumption, we have $\md(a_i\vee_La_j)=\md(\Delta)=a_0$. Thus we get the result on $\Xx_2$. Lastly, if $[\alpha_1,\cdots,\alpha_n]$ is an $n$-cell for $n\geqslant 2$, then $[\alpha_{n-1},\alpha_{n}]$ is a $2$-cell. Thus $\alpha_{n-1}=a_0$, and $\alpha_{n-2}<a_0$ is impossible if $n>2$. We obtain that $\Xx_n=\varnothing$ if $n\geqslant 3$.\end{proof}

Let $M=M(m,\ell)$ be a circular monoid. We know that $H_0(G(m,\ell),\Z)=\Z$ and that $H_1(G(m,\ell),\Z)$ is the abelianization $G(m,\ell)^{\ab}$ of $G(m,\ell)$. Since $C_n=\{0\}$ for $n\geqslant 3$, we have $H_n(G(m,\ell),\Z)=\{0\}$ for $n\geqslant 3$. Furthermore, the group $H_2(G(m,\ell),\Z)$ is the kernel of the map $\partial_2:C_2\otimes \Z\to C_1\otimes \Z$. In particular, it is a free abelian group. Since the Euler characteristic of the complex $(C_n,\partial_n)$ is $0$, we get that $H_2(G(m,\ell),\Z)\simeq \Z^{r-1}$ where $r$ is the rank of the free part of $H_1(G(m,\ell),\Z)$. Thus the integral homology of $G(m,\ell)$ can be computed by only computing $G(m,\ell)^{\ab}$.

\begin{lem}\label{lem:homology_circular}
Let $m,\ell$ be positive integers. We have $G(m,\ell)^{\ab}\simeq \Z^{m\wedge \ell}$. The integral homology of $G(m,\ell)$ is then given by
\[H_n(G(m,\ell),\Z)=\begin{cases} \Z&\text{if }n=0\\ \Z^{m\wedge \ell} &\text{if }n=1\\ \Z^{m\wedge \ell-1}&\text{if }n=2\\ 0&\text{if }n\geqslant 3\end{cases}\]
\end{lem}
\begin{proof}
Let $a_0,\ldots,a_{m-1}$ denote the atoms of $M(m,\ell)$. In $G(m,\ell)^{\ab}$, we have $a_{i+\ell}=s(i,\ell)^{-1}a_is(i,\ell)=a_i$ for all $i\in \intv{0,m-1}$. Conversely, in the group $\Z^{m}$ quotiented by the relations $a_i=a_{i+\ell}$ for all $i\in \intv{0,m-1}$, we have
\[s(i+1,\ell)=s(i+1,\ell-1)a_{i+\ell}=s(i+1,\ell-1)a_i=a_is(i+1,\ell-1)=s(i,\ell)\]
for all $i\in \intv{0,m-1}$. Thus we have
\[G(m,\ell)^{\ab}=\left\langle a_0,\ldots,a_{m-1}~\left|~ \begin{cases} a_i=a_{i+\ell}&\forall i\in \intv{0,m-1}\\ a_ia_j=a_ja_i&\forall i,j\in \intv{0,m-1} \end{cases}\right.\right\rangle.\]
This group is free abelian, with rank the cardinality of $(\Z/m\Z)/(\ell\Z/m\Z)=\Z/(m\wedge \ell)\Z$.\end{proof}
In particular, for $m=2$, we recover the result of \cite[Table 1]{salvetti} on the homology of spherical Artin groups of rank $2$.

\subsubsection{Remarkable isomorphisms.}\label{sec:isomorphisms}
In this section we give the classification of circular groups up to group isomorphism. First, Corollary \ref{cor:center_circular_group} gives that, for any positive integers $m,\ell$, the group $G(m,\ell)$ is abelian if and only if $m=1$ or $\ell=1$ or $m=\ell=2$, in which case $G(m,\ell)\simeq \Z$, $G(m,\ell)\simeq \Z$ and $G(m,\ell)\simeq \Z^2$, respectively.

\begin{prop}\label{prop:isomorphism_restriction_circular}
Let $m,\ell,m',\ell'$ be four positive integers. If the groups $G(m,\ell)$ and $G(m',\ell')$ are isomorphic and nonabelian, then $(m',\ell')\in \{(m,\ell),(\ell,m)\}$.
\end{prop} 
\begin{proof}
Let $d:=m\wedge \ell$ and $d':=m'\wedge \ell'$. If $G(m,\ell)$ and $G(m',\ell')$ are isomorphic, then Lemma \ref{lem:homology_circular} gives that
\[\Z^d\simeq H_1(G(m,\ell),\Z)\simeq H_1(G(m',\ell'),\Z)\simeq \Z^{d'}.\]
In particular we have $d'=d$.

Since $G(m,\ell)$ and $G(m',\ell')$ are nonabelian, both the centers of $G(m,\ell)$ and $G(m',\ell')$ are cyclic and generated by some power of $\Delta$. An isomorphism $f:G(m,\ell)\to G(m',\ell')$ then induces a bijection between irreducible periodic elements of $G(m,\ell)$ and of $G(m',\ell')$. By Proposition \ref{prop:irreducible_periodic_circulaire}, we have
\begin{enumerate}[-]
\item If $m|\ell$, then we either have $m'|\ell'$, in which case we have $\frac{\ell}{m}=\frac{\ell'}{m'}$, $d=m=m'=d'$ and $\ell=\ell'$, or $\ell'|m'$, in which case we have $\frac{\ell}{m}=\frac{m'}{\ell'}$, $d=m=\ell'=d'$ and $\ell=m'$.
\item If $\ell|m$, the same reasoning gives $(m,\ell)=(m',\ell')$ or $(m,\ell)=(\ell',m')$.
\item Lastly, if neither $m|\ell$ nor $\ell|m$, then we have neither $m'|\ell'$ nor $\ell'|m'$. We then have either
\[\frac{\ell}{d}=\frac{\ell'}{d'} \text{ and }\frac{m}{d}=\frac{m'}{d'}\text{ or } \frac{\ell}{d}=\frac{m'}{d'} \text{ and }\frac{m}{d}=\frac{\ell'}{d'}.\]
Since $d'=d$, we obtain $(m',\ell')\in \{(m,\ell),(\ell,m)\}$.
\end{enumerate}
\end{proof}

This proposition strongly restricts the possible isomorphisms between circular groups. We can then show that all the remaining possible isomorphisms actually occur:

\begin{prop}\label{prop:isomorphism_circular_group}
Let $m,\ell$ be two positive integers. There is an isomorphism of groups between $G(m,\ell)$ and $G(\ell,m)$, which sends atoms of $M(m,\ell)$ to conjugates of atoms in $M(\ell,m)$.
\end{prop}
\begin{proof}
The result is immediate if $m=\ell$. Up to exchanging $m$ and $\ell$, we can assume that $m<\ell$. Let $\ell=mp+r$ be the Euclidean division of $\ell$ by $m$. 

We denote by $\{a_0,\ldots,a_{m-1}\}$ the atoms of $M(m,\ell)$ and by $\{b_0,\ldots,b_{\ell-1}\}$ the atoms of $M(\ell,m)$. We also consider $F_m$ to be the free group generated by $\{a_0,\ldots,a_{m-1}\}$. Exceptionally, we denote the simple elements of $M(\ell,m)$  by $t(i,p)$ instead of $s(i,p)$ to avoid confusions with the simple elements of $M(m,\ell)$. We also denote by $\ttilde{s}(i,p)$ the product $a_i\cdots a_{i+p}$ in $F_m$.

Let $\ttilde{f}:F_m\to G(\ell,m)$ be the morphism defined by
\[\begin{cases} f(a_0):=b_{m-1},\\ f(a_i):=(b_{m-i-1})^{f(\ttilde{s}(0,i))}&\forall i\in \intv{1,m-1}.\end{cases}\]
By an immediate induction, we get that $\ttilde{f}(\ttilde{s}(0,k))=t(m-k,k)$ for all $k\in \intv{0,m-1}$. We show that $\ttilde{f}$ induces a well defined group morphism $f:G(m,\ell)\to G(\ell,m)$. Let $i\in \intv{1,m-1}$, we have
\begin{align*}
\ttilde{f}(\ttilde{s}(i,\ell))&=\ttilde{f}(\ttilde{s}(0,i)^{-1}\ttilde{s}(0,\ell+i))\\
&=\ttilde{f}(\ttilde{s}(0,i)^{-1})\ttilde{f}(\ttilde{s}(0,\ell)) \ttilde{f}(\ttilde{s}(0,r+i))\\
&=t(m-i,i)^{-1}\ttilde{f}(\ttilde{s}(0,m)^q)) t(m-r-i,r+i)\\
&=t(m-i,i)^{-1} t(0,m)^qt(m-r-i,r+i)\\
&=t(0,m)^q t(m-i+qm,i)^{-1}t(m-r-i,r+i)\\
&=t(0,m)^q t(m-i-r,i)^{-1}t(m-r-i,r+i)\\
&=t(0,m)^q t(0,r)=\ttilde{f}(\ttilde{s}(0,\ell)).
\end{align*}
We show that $f$ is an isomorphism by constructing its inverse. First, by definition of $f$, we have
\[\forall i\in \intv{0,m-1},~b_{m-i-1}=f({}^{s(0,i)}a_i)\]
and we define $g(b_j)=a_{m-j-1}^{s(0,m-j-1)^{-1}}$ for $j\in \intv{0,m-1}$. We also have $\Delta=t(0,m)=f(s(0,m))$ and we define $g(t(0,m))=s(0,m)$. Let $j\in \intv{0,\ell-1}$ and let $j=mp+j'$ be the Euclidean division of $j$ by $m$. We have $b_j=\phi^p(b_{j'})=\Delta^{-p}b_{j'}\Delta^p$ and we define
\[g(b_j)=g(\Delta^p b_{j'}\Delta^p):=s(0,m)^{-p}g(b_{j'})s(0,m)^p=a_{m-j'-1}^{s(m-j'-1,j'+1)s(0,m)^{p-1}}\]
To show that $g$ does define a group morphism $G(\ell,m)\to G(m,\ell)$, we have to show that $g(t(i,m))$ does not depend on $i$. We have $g(t(0,m))=s(0,m)$ by definition. Then, let $i\in \intv{1,\ell-1}$ be such that $g(t(i-1,m))=s(0,m)$. Let $i-1=mp+k$ be the Euclidean division of $i-1$ by $m$. We have
\begin{align*}
g(t(i,m))&=g(b_{i-1}^{-1})g(t(i-1,m))g(b_{i+m-1})=s(0,m)\\
&=a_{m-k-1}^{-1~s(m-k-1,k+1)s(0,m)^{p-1}} s(0,m) a_{m-k-1}^{s(m-k-1,k+1)s(0,m)^{p}}\\
&=s(0,m).
\end{align*}
We obtain that $g(t(i,m))=g(t(0,m))=s(0,m)$ by induction. It is an immediate check to see that $f$ and $g$ are inverse to each other.
\end{proof}

If we combine Proposition \ref{prop:isomorphism_restriction_circular} and Proposition \ref{prop:isomorphism_circular_group}, we get a complete classification of circular groups up to group isomorphism.
\begin{cor}\label{cor:classification_circular_groups}
Let $m,\ell,m',\ell'$ be four positive integers. The groups $G(m,\ell)$ and $G(m',\ell')$ are isomorphic if and only if one of the following holds
\begin{enumerate}[-]
\item $1\in \{m,\ell\}$ and $1\in \{m',\ell'\}$. In this case, $G(m,\ell)\simeq G(m',\ell')\simeq \Z$.
\item $(m',\ell')\in \{(m,\ell),(\ell,m)\}$.
\end{enumerate}
\end{cor}

\begin{exemple}
If $m=2$, then $M(2,\ell)$ is the Artin monoid for the Artin group of type $I_2(\ell)$, while $M(\ell,2)$ is the dual braid monoid for the same Artin group. The isomorphism $G(2,\ell)\to G(\ell,2)$ constructed in the above proof is already known: it sends $a_0$ to $b_1$ and $a_1$ to $b_0^{b_1}=b_1^{-1}\Delta=b_2$.
\end{exemple}

\subsection{Application to complex braid groups of rank 2} We refer the reader to \cite{lehrertaylor} for general results on complex reflection groups. They are finite subgroups of $\GL_n(\C)$ generated by complex (pseudo)-reflections. 

The \nit{braid group} associated to a complex reflection group $W$ is defined as $\pi_1(X/W)$, where $X$ is the complement inside $\C^n$ of the hyperplane arrangement defined by the reflections of $W$. The classification of irreducible complex reflection groups was done in \cite{shetod}. It separates irreducible complex reflection groups into a general series $G(de,e,n)$ depending on integer parameters $d,e,n$ and a list of $34$ exceptional cases $G_{4},\ldots,G_{37}$.

\begin{fact*}
Let $W\subset \GL_2(\C)$ be a complex reflection group of rank $2$. The braid group $B(W)$ is isomorphic to a circular group.
\end{fact*}
This is mostly a rephrasing of \cite[Theorem 1 and Theorem 2]{bannai}: 
\begin{enumerate}[-]
\item The only non-irreducible cases are groups of the form $W=\Z/d\Z/\times \Z/d'\Z$ (with $d,d'\geqslant 1$). In this case we have $B(W)\simeq \Z\times \Z\simeq G(2,2)$.
\item If $W=G(de,e,2)$ for $e$ odd and $d\geqslant 2$ or $W\in \{G_5,G_{10},G_{18}\}$, then $B(W)\simeq G(2,4)$.
\item If $W=G(de,e,2)$ for $e$ even and $d\geqslant 2$ or $W\in \{G_7,G_{11},G_{15},G_{19}\}$, then $B(W)\simeq G(3,3)$.
\item If $W=G(e,e,2)$ for $e\geqslant 3$, then $B(W)\simeq G(2,e)$.
\item If $W\in \{G_4,G_8,G_{16}\}$, then $B(W)\simeq G(2,3)$.
\item If $W\in \{G_6,G_9,G_{13},G_{17}\}$, then $B(W)\simeq G(2,6)$.
\item If $W=G_{14}$, then $B(W)\simeq G(2,8)$.
\item If $W=G_{20}$, then $B(W)\simeq G(2,5)$.
\item If $W=G_{21}$, then $B(W)\simeq G(2,10)$.
\item If $W=G_{12}$, then $B(W)\simeq G(3,4)$.
\item If $W=G_{22}$, then $B(W)\simeq G(3,5)$.
\end{enumerate}

A direct application of Theorem \ref{theo:roots_conjugate_circular} then gives
\begin{theo}\label{theo:roots_conjugate_rank2}
Let $W$ be a complex reflection group of rank $2$, and let $B(W)$ be its braid group. If $\alpha,\beta\in B(W)$ are such that $\alpha^n=\beta^n$ for some nonzero integer $n$, then $\alpha$ and $\beta$ are conjugate in $B(W)$.
\end{theo}

\begin{rem}
Our approach only covers complex reflection groups of rank $2$. Indeed, by Lemma \ref{lem:cellules_dehlaf_circular}, circular groups have homological dimension at most 2, and a complex braid group of rank $r$ has homological dimension $r$ by \cite[Proposition 1.1]{homcomp1}. Thus only rank 2 complex reflection groups can have a braid group isomorphic to a circular group.
\end{rem}

\section{$\Delta$-product and hosohedral-type groups}
In this section we present hosohedral-type Garside groups as a generalization of circular groups. These groups are enveloping groups of so-called hosohedral-type monoids. These monoids were first introduced by Picantin in his PhD thesis under the name ``monoïdes de type fuseau'' (\cite[Définition 1.3]{thespicantin}). In \cite[Proposition 2.4]{thespicantin}, Picantin shows that these monoids are exactly the Garside monoids whose lattice of simples has the shape of a hosohedron (``fuseau'' in french). The name ``hosohedral-type monoid'' was suggested to us by Picantin. 

More recently, hosohedral-type groups were identified by Mireille Soergel in \cite[Theorem 4.6]{soergelsystolic} as the Garside groups satisfying a particular nonpositive curvature property (namely, the systolicity of the flag complex associated to the ``Garside presentation'' as in \cite[Lemma 4.3]{soergelsystolic}).

Here we introduce these groups as a particular case of a general construction, already present in \cite{dehpar}, which we call the $\Delta$-product of Garside monoids.

\subsection{$\Delta$-product of Garside monoids}
Let $(M_1,\Delta_1),\ldots,(M_h,\Delta_h)$ be a finite family of homogeneous Garside monoids that we fix throughout this section. We denote by $\ell_1,\ldots,\ell_h$ the associated length functions, and by $S_1,\ldots,S_h$ the associated set of simples. We also set $A_1,\ldots,A_h$ the set of atoms of $M_1,\ldots,M_h$, respectively.

The free product $M_1*\cdots *M_h$ is not a Garside monoid because two atoms coming from a different factor $M_i$ do not have any common multiple (in particular, no lcms). We can fix this by forcing the $\Delta_i$ to all be equal.

\begin{definition}
Let $(M_1,\Delta_1),\ldots,(M_h,\Delta_h)$ be a family of homogeneous Garside monoids. The \nit{$\Delta$-product} of the $M_i$ is defined by 
\[M_1*_\Delta M_2*_\Delta\cdots *_\Delta M_h:=\bigslant{*_{i=1}^h M_i}{(\Delta_i=\Delta_j~\forall i,j\in \intv{1,h})}\]
Likewise, we define the $\Delta$-product of the enveloping groups $G(M_i)$ by
\[G(M_1)*_\Delta G(M_2)*_\Delta\cdots *_\Delta G(M_h)=\bigslant{*_{i=1}^h G(M_i)}{(\Delta_i=\Delta_j~\forall i,j\in \intv{1,h})}\]
\end{definition}

\begin{rem}
The definition of the $\Delta$-product of Garside monoids depends really on the Garside element and not only on the monoid themselves. For instance, we have $M(1,\ell)\simeq \Z_{\geqslant 0}$ as a monoid for every $\ell\geqslant 1$. However, we have
\[M(1,p)*_\Delta M(1,q)\simeq \langle a,b~|~a^p=b^q\rangle^+.\]
If $p,q\geqslant 2$, then this monoid has two atoms and cannot be isomorphic to $\Z_{\geqslant 0}\simeq M(1,1)*_\Delta M(1,1)$. Furthermore, the enveloping group of this monoid is not always isomorphic to $\Z=G(M(1,1)*_\Delta M(1,1))$.
\end{rem}

\begin{rem}\label{rem:identity_delta_product}
Note that, if $(M,\Delta)$ is a Garside monoid, the monoid $M*_\Delta M(1,1)$ is naturally isomorphic to $M$. Thus we can assume that all the $(M_i,\Delta_i)$ are distinct from $M(1,1)$.
\end{rem}

We first show that the enveloping group of the $\Delta$-product of the $(M_i,\Delta_i)$ identifies with the $\Delta$-product of the enveloping groups $G(M_i)$.

\begin{lem}
Let $(M_1,\Delta_1),\ldots,(M_h,\Delta_h)$ be a family of homogeneous Garside monoids. Let also $M$ be the $\Delta$-product $M_1*_\Delta\cdots*_\Delta M_h$. The groups $G(M)$ and $G(M_1)*_\Delta \cdots *_\Delta G(M_h)$ are naturally isomorphic.
\end{lem}
\begin{proof}
If $M,M'$ are two monoids, we denote by $\ho(M,M')$ the set of monoid morphisms from $M$ to $M'$. If $G$ and $G'$ are two groups, $\ho(G,G')$ is in fact the set of group morphisms from $G$ to $G'$. Let $H$ be a group. By definition of the enveloping group and of the $\Delta$-product, we have natural bijections
\begin{align*}
\ho(G(M),H)&\simeq \ho(M_1*_\Delta\cdots *_\Delta M_h,H)\\
&\simeq \{f\in \ho(M_1*\cdots*M_h,H)~|~ \forall i,j\in \intv{1,h}, f(\Delta_i)=f(\Delta_j)\}\\
&\simeq \left\{(f_i) \in \prod_{i=1}^n \ho(M_i,H)~\left|~ \forall i,j, f_i(\Delta_i)=f_j(\Delta_j)\right. \right\}\\
&\simeq \left\{(f_i) \in \prod_{i=1}^n \ho(G(M_i),H)~\left|~ \forall i,j, f_i(\Delta_i)=f_j(\Delta_j)\right. \right\}\\
&\simeq \{f\in \ho(G(M_1)*\cdots*G(M_h),H)~|~ \forall i,j\in \intv{1,h}, f(\Delta_i)=f(\Delta_j)\}\\
&\simeq \ho(G(M_1)*_\Delta\cdots *_\Delta G(M_h),H)
\end{align*}
Applying this to $H:=G(M)$ gives a bijection $\ho(G(M),G(M))\simeq \ho(G(M_1)*_\Delta\cdots *_\Delta G(M_h),G(M))$. The image of the identity morphism of $G(M)$ under this bijection gives the desired result.
\end{proof}

From now on, we will identify $G(M_1*_\Delta\cdots*_\Delta M_h)$ with $G(M_1)*_\Delta\cdots *_\Delta G(M_h)$. We now fix a family $(M_1,\Delta_1),\ldots,(M_h,\Delta_h)$ of Garside monoids, all distinct from $M(1,1)$. We denote by $M$ the associated $\Delta$-product, and by $A$ its set of atoms. Let $i\in \intv{1,h}$. By definition of the $\Delta$-product as a quotient of the free product, there is a natural morphism $\pphi_i:G(M_i)\to G(M)$. We also denote by $\pphi_i$ the restriction from $M_i$ to $M$.

\begin{prop}\cite[Proposition 5.3]{dehpar}\label{prop:simples_delta_product}\newline
The atoms of the monoid $M$ are given by $A:=\bigsqcup_{i=1}^h \pphi_i(A_i)$. Furthermore, $(M,\Delta)$ is a homogeneous Garside monoid with simple elements
\[S=\bigcup_{i=1}^h \pphi_i(S_i).\]
Where $\Delta=\pphi_i(\Delta_i)$ for any $i\in \intv{1,h}$. The Garside automorphism $\phi$ of $(M,\Delta)$ is given on a simple $\pphi_i(s)$ by $\phi(\pphi_i(s))=\pphi_i(\phi_i(s))$, where $\phi_i$ is the Garside automorphism of $(M_i,\Delta_i)$.
\end{prop}
The assertion on the atoms and the assertion that $(M,\Delta)$ is a homogeneous Garside monoid are in the original statement of \cite[Proposition 5.3]{dehpar}. For $i\in \intv{1,h}$, we identify $A_i$ with the subset $\pphi_i(A_i)$ of $A$ from now on. The proof of \cite[Proposition 5.3]{dehpar} can also be used to show the assertion on the simple elements:

First, let $i\in \intv{1,h}$ and let $s\in S_i$ be a simple element of $(M_i,\Delta_i)$. We have $s\bbar{s}=\Delta_i$ and $\pphi_i(s)\pphi_i(\bbar{s})=\Delta$. Thus $\pphi_i(s)$ is simple and $\bbar{\pphi_i(s)}=\pphi_i(\bbar{s})$. We also obtain that the Garside automorphism is given by $\phi(\pphi_i(s))=\bbar{\bbar{\pphi_i(s)}}=\pphi_i(\phi_i(s))$. Conversely, we have to show that, if $s\in S$ is a simple element of $(M,\Delta)$, then there is some $i\in \intv{1,h}$ and some simple $\ttilde{s}\in S_i$ with $s=\pphi_i(\ttilde{s})$. This is a direct consequence of the following lemma.

\begin{lem}\label{lem:injectivity_of_pphi_i} Let $w$ be a word in $A$ which expresses a simple element $s\in S$. There is some $i\in \intv{1,h}$ such that all the letters of $w$ lie inside $\pphi_i(A_i)$. The word $w$ then also represents some $\ttilde{s}\in S_i$ with $\pphi_i(\ttilde{s})=s$. Furthermore, if $s\notin\{1,\Delta\}$, then the integer $i$ and the simple $\ttilde{s}\in S_i$ are unique.
\end{lem}

\begin{proof}
Let $i\in \intv{1,h}$. For each pair $a,b\in A_i$, we choose two words $f_i(a,b)$ and $f_i(b,a)$ such that $af_i(a,b)$ and $bf_i(b,a)$ are two words in $A_i$ expressing $a\vee b$ in $M_i$. By \cite[Theorem 4.1]{dehpar}, the monoid $M_i$ admits the following presentation
\[M_i=\langle af_i(a,b)=bf_i(b,a)~|~ a,b\in A_i\rangle^+.\]
Now, for $a\in A_i$, we choose a word $c(a)$ in $A_i$ representing $\bbar{a}$ in $M_i$. The proof of \cite[Proposition 5.3]{dehpar} gives that $M$ admits the presentation
\begin{equation}\label{eq:pres_delta_prod}M\simeq \left\langle af(a,b)=bf(b,a)~\left| f(a,b)=\begin{cases} f_i(a,b)&\text{if }a,b\in A_i\\ c(a)&\text{if }a\in A_i,b\in A_j,i\neq j\end{cases}\right. \right\rangle^+.\end{equation}

By definition, if $w_1=w_2$ is a relation in this presentation, then we either have
\begin{enumerate}[-]
\item There is an $i\in \intv{1,h}$ such that all the letters of both $w_1$ and $w_2$ lie in $A_i$. In this case, $w_1=w_2$ also holds in $M_i$.
\item There are two distinct integers $i,j\in \intv{1,h}$ and two atoms $a\in A_i,b\in A_j$ such that $w_1=ac(a)$ and $w_2=bc(b)$. In this case, $w_1$ (resp. $w_2$) represents $\Delta_i$ in $M_i$ (resp. $\Delta_j$ in $M_j$).
\end{enumerate}

Let $w$ be a word in $A$ which represents $\Delta$ in $M$, and let $a\in A$ be an atom of $M$. By assumption, there is a sequence of words $w_1,\ldots,w_m$ in $A$ such that $w_1=w$, $w_m=ac(a)$ and each $w_k$ is equivalent to $w_{k+1}$ by the use of one relation of presentation (\ref{eq:pres_delta_prod}) for $k\in \intv{1,m-1}$. Up to changing the atom $a$, we can assume that $m$ is the first integer such that $w_m$ is equal to a word of the form $bc(b)$ for some atom $b$. 

Let $i\in \intv{1,h}$ be such that $a\in A_i$. We claim that for all $j\in \intv{1,m}$, the word $w_j$ contains letters only in $A_i$. If this is not the case, then let $k_0$ be the last integer in $\intv{1,m}$ such that $w_{k_0}$ contains letters not lying in $\pphi_i(A_i)$ (we have $k_0<m$ by assumption). The defining relations of (\ref{eq:pres_delta_prod}) giving $w_{k_0}=w_{k_0+1}$ in $M$ then have the form $bc(b)=ac(a)$. Since $w_{k_0}$ and $bc(b)$ both express $\Delta$ in $M$, we have $w_{k_0}=bc(b)$, which contradicts the minimality assumption on $m$. 

Since all the letters of $w_j$ lie inside $A_i$ for $j\in \intv{1,m}$, we have that the relations of (\ref{eq:pres_delta_prod}) giving $w_j=w_{j+1}$ for $j\in \intv{1,m-1}$ also hold in $M_i$. Thus $w_1$ is a word in $A_i$, which expresses the element $\Delta_i$ in $M_i$.

Let now $s\in S$ be a simple element of $M$. By definition, there is a word $w_1=ww_2$ in $A$ expressing $\Delta$ such that $w$ expresses $s$ in $M$. The first part of the proof gives that there is some $i\in \intv{1,h}$ such that $w_1,w$ and $w_2$ are actually words in $A_i$. The word $w$ (resp. $w_2$) then expresses an element $\ttilde{s}$ (resp. $\ttilde{s}'$) in $M_i$, with $\ttilde{s}\ttilde{s}'=\Delta_i$. We then have $\ttilde{s}\in S_i$ and $\pphi_i(\ttilde{s})=s$.

Suppose now that $s\notin\{1,\Delta\}$, and let $w'$ be another word in $A$ expressing $s$ in $M$. By definition, there is a sequence of words $w_1,\ldots,w_n$ of words in $A$, such that $w_1=w$, $w_n=w'$ and each $w_k$ is equivalent to $w_{k+1}$ by the use of one relation of presentation (\ref{eq:pres_delta_prod}). Since $s\prec \Delta$, none of the $w_k$ contain a subword expressing $\Delta$ in $M$. We deduce that all the relations giving $w_k=w_{k+1}$ in $M$ for $k\in \intv{1,n-1}$ are between words in $A_i$ and also hold in $M_i$. Thus $w_k$ also expresses $\ttilde{s}$ for all $k\in \intv{1,n}$.
%
%If $w'$ contains letters not lying in $A_i$, then let $k$ be the smallest $k$ such that $w_{k+1}$ contains letters not lying in $A_i$. In this case, the defining relation of (\ref{eq:pres_delta_prod}) giving that $w_k=w_{k+1}$ in $M$ has the form $ac(a)=bc(b)$. We then obtain that $s$ can be written as $t\Delta u$ in $M$. Since $s$ is a simple elements, this implies that $s=\Delta$, which is a contradiction. 
%
%
%
%
%We deduce that $w_k$ is a word in $A_i$ for all $k\in \intv{1,n}$, in particular $w=w_n$ is a word in $A_i$. Furthermore, we also obtain that all the relations giving $w_k=w_{k+1}$ in $M$ for $k\in \intv{1,n-1}$ also hold in $M_i$, thus $w_k$ also expresses $\ttilde{s}$ for all $k\in \intv{1,n}$.
\end{proof}

\begin{prop}\label{prop:4.9}Let $i\in \intv{1,n}$ and let $s\in S_i$. The morphism $\pphi_i$ induces an injective morphism of lattices from $\{t\in S_i~|~ t\preceq s\}$ to $\{t\in S~|~ t\preceq \pphi_i(s)\}$. Furthermore, if $s\neq \Delta_i$, then this morphism is also surjective.
\end{prop}
\begin{proof}
We show the injectivity and surjectivity assumptions before considering lcms and gcds. Let $s\in S_i$ and let $X$ and $Y$ denote the two considered sets. Since $\pphi_i$ is a morphism of monoids, it restricts to a poset morphism $X\to Y$. First, we show that $\pphi_i:X\to Y$ is always injective. Let $t,t'\in X$ be such that $\pphi_i(t)=\pphi_i(t')$. Let $w$ and $w'$ be two words in $A_i$ expressing $t$ and $t'$, respectively. The words $w$ and $w'$ express the same element $\pphi_i(t)$ in $M$. If $\pphi_i(t)=\Delta$, then $w$ and $w'$ are two words in $A_i$ representing $\Delta$ in $M$. The proof of Lemma \ref{lem:injectivity_of_pphi_i} gives that $w$ and $w'$ represent the same element in $M_i$, which is $\Delta_i$. If $\pphi_i(t)\neq \Delta$, then $w$ and $w'$ express the same element in $S_i$ by Lemma \ref{lem:injectivity_of_pphi_i}, thus $t=t'$.

Now, suppose that $s\neq \Delta$. We show that $\pphi_i:X\to Y$ is surjective. Let $t\preceq \pphi_i(s)$ in $M$. By definition, there is a word $w=w_1w_2$ in $A$ expressing $\pphi_i(s)$ such that $w_1$ expresses $t$. By Lemma \ref{lem:injectivity_of_pphi_i}, $w$ and $w_1$ are words in $A_i$. The word $w_1$ then expresses some element $\ttilde{t}\in S_i$ such that $\pphi_i(\ttilde{t})=t$ and $\pphi_i:X\to Y$ is surjective.

We now show that $\pphi_i:S_i\to S$ is a morphism of lattices. That is, $\pphi_i$ preserves right-lcms and left-gcds. Let $s,t\in S_i$ be two simples of $(M_i,\Delta_i)$, and let $u=s\wedge t$. Let $x:=\pphi_i(s)\wedge \pphi_i(t)$. Since the simple $\pphi_i(u)$ is obviously a left-divisor of $\pphi_i(s)$ and $\pphi_i(t)$, we have $\pphi_i(u)\preceq x$. If $s=\Delta_i$, then $u=t$ and $x=\pphi_i(t)=\pphi_i(u)$. Likewise the result is clear if $t=\Delta_i$. We assume from now on that $t,s\neq \Delta_i$. Since $s\neq \Delta_i$ (resp. $t\neq \Delta_i$), the first part of the proof gives the existence of a unique $\ttilde{x}$ (resp. $\ttilde{x}'$) in $S_i$ such that $\ttilde{x}\preceq s$ (resp. $\ttilde{x}'\preceq t$) and $\pphi_i(\ttilde{x})=x=\pphi_i(\ttilde{x}')$. Since $\pphi_i$ is injective on $S_i$, we get that $\ttilde{x}=\ttilde{x}'$ is a common divisor of $s$ and $t$ in $M_i$. We obtain $\ttilde{x}\preceq u$, $x\preceq \pphi_i(u)$ and $x=\pphi_i(u)$.

Let now $v=s\vee t$ and $y=\pphi_i(s)\vee \pphi_i(t)$. Again, since $\pphi_i(v)$ is an obvious right-multiple of both $s$ and $t$, we have $y\preceq \pphi_i(v)$. By the first part of the proof, there is a unique $\ttilde{y}\in S_i$ such that $\ttilde{y}\preceq v$ and $\pphi_i(\ttilde{y})=y$. The first part of the proof also gives that $s,t\preceq \ttilde{y}$. Thus $v\preceq \ttilde{y}$, $\pphi_i(v)\preceq \pphi_i(\ttilde{y})=y$ and $\pphi_i(v)=y$.
\end{proof}

Of course, one can show by similar arguments that $\pphi_i$ is an injective morphism of lattices from $(S_i,\succeq)$ to $(S,\succeq)$.

\begin{cor}\label{cor:lcm_gcd_delta_product}
Let $s,t$ be two simple elements in $M$, both different from $\Delta$ and $1$. Assume that $s=\pphi_i(\ttilde{s})$ and $t=\pphi_j(\ttilde{t})$ for $i\neq j$, $\ttilde{s}\in S_i$ and $\ttilde{t}\in S_j$. We have $s\wedge t=1$ and $s\vee t=\Delta$ in $M$. Furthermore, the word $st$ is greedy in $(M,\Delta)$.
\end{cor}
\begin{proof} 
Let $u:=s\wedge t$. By Proposition \ref{prop:4.9}, there is a unique $\ttilde{u}\in S_i$ (resp. $\ttilde{u}'\in S_j$) such that $\pphi_i(\ttilde{u})=\pphi_j(\ttilde{u}')=u$. Since $i\neq j$, Lemma \ref{lem:injectivity_of_pphi_i} gives that $u=1$. We apply similar reasoning for lcms. This gives in particular that $\bbar{s}\wedge t=1$, thus the path $st$ is greedy.
\end{proof}

\begin{prop}\label{prop:left-weighted_facto_in_delta_product}
Let $i\in \intv{1,h}$. The morphism $\pphi_i:G(M_i)\to G(M)$ preserves left-weighted factorizations. In particular it is injective.
\end{prop}
\begin{proof}
We first show that $\pphi_i:M_i\to M$ preserves greedy normal forms. By definition, we only have to show that, if $st$ is a greedy word in $M_i$, then $\pphi_i(s)\pphi_i(t)$ is a greedy word in $M$. If $st$ is a greedy word of length $2$ in $M_i$, then Proposition \ref{prop:4.9} gives
\[\bbar{\pphi_i(s)}\wedge \pphi_i(t)=\pphi_i(\bbar{s})\wedge \pphi_i(t)=\pphi_i(\bbar{s}\wedge t)=1.\]
The word $\pphi_i(s)\pphi_i(t)$ is then greedy by Lemma \ref{gr2}.

Let now $x=\Delta_i^ks_1\cdots s_r$ be the left-weighted factorization of some $x\in M_i$. Since $\Delta_i\not\preceq s_1\cdots s_r$, we have $\Delta\not\preceq \pphi_i(s_1\cdots s_r)$. Furthermore, the word $\pphi_i(s_1)\cdots \pphi_i(s_r)$ is greedy because $\pphi_i:M_i\to M$ preserves greediness. The word $\Delta^k \pphi_i(s_1)\cdots \pphi_i(s_r)$ is then the left-weighted factorization of $\pphi_i(x)$ by definition.
\end{proof}

\begin{prop}\label{prop:sss_in_delta_product}
Let $i\in \intv{1,h}$, and let $x\in G(M_i)$. We have $\pphi_i(\SSS(x))=\SSS(\pphi_i(x))$. Furthermore, if $x$ is not conjugate to a power of $\Delta_i$ in $M_i$, then the centralizer of $\pphi_i(x)$ in $G(M)$ is the image under $\pphi_i$ of the centralizer of $x$ in $G(M_i)$.
\end{prop}
\begin{proof}
Let $x$ be an element of $G(M_i)$. By Proposition \ref{prop:left-weighted_facto_in_delta_product}, we have $\pphi_i(\init(x))=\init(\pphi_i(x))$ and $\pphi_i(\fin(x))=\fin(\pphi_i(x))$. Thus, we also have $\pphi_i(\cyc(x))=\cyc(\pphi_i(x))$ and $\pphi_i(\dec(x))=\dec(\pphi_i(x))$. 

If $x$ is conjugate to some $\Delta_i^k$ with $k\in \Z$, then we have $\SSS(x)=\{\Delta^k\}=\SSS(\pphi_i(x))$. 

From now on, we suppose that $x$ is not conjugate to any power of $\Delta$. We show that $\inf(\SSS(x))=\inf(\SSS(\pphi_i(x))$ and $\sup(\SSS(x))=\sup(\SSS(\pphi_i(x))$. First, let $x'\in \SSS(x)$, Proposition \ref{prop:left-weighted_facto_in_delta_product} gives that
\[\inf(\SSS(\pphi_i(x)))\geqslant\inf(\pphi_i(x'))=\inf(x')=\inf(\SSS(x)),\]
\[\sup(\SSS(\pphi_i(x)))\leqslant \sup(\pphi_i(x'))=\sup(x')=\sup(\SSS(x)).\]
Conversely, one can reach an element $y$ of $\SSS(\pphi_i(x))$ by applying a sequence of cycling and decycling to $\pphi_i(x)$. Applying the same operations to $x$ gives a conjugate $\ttilde{y}$ of $x$ in $G(M_i)$ such that $\pphi_i(\ttilde{y})=y$. We then have
\[\inf(\SSS(x))\geqslant \inf(\ttilde{y})=\inf(y)=\inf(\SSS(\pphi_i(x))),\]
\[\sup(\SSS(x))\leqslant \sup(\ttilde{y})=\sup(y)=\sup(\SSS(\pphi_i(y))).\]
Thus, $\inf(\SSS(x))=\inf(\SSS(\pphi_i(x)))$ and $\sup(\SSS(x))=\sup(\SSS(\pphi_i(x)))$ as claimed.

Now, we show that $\pphi_i(\SSS(x))\subset \SSS(\pphi_i(x))$. Let $x'\in \SSS(x)$, we have $\inf(\pphi_i(x'))=\inf(x')=\inf(\SSS(x))=\inf(\SSS(\pphi_i(x)))$, and likewise, $\sup(\pphi_i(x'))=\sup(\SSS(\pphi_i(x)))$. Since $\pphi_i(x')$ is conjugate to $\pphi_i(x)$, we have $\pphi_i(x')\in \SSS(\pphi_i(x))$.

Let $x'\in \SSS(x)$, and let $s\in S$ be a simple element of $M$ such that $\pphi_i(x')^s\in \SSS(\pphi_i(x))$. We show that $s\in \pphi_i(S_i)$. Let $x'=\Delta_i^k s_1\cdots s_r$ be the left-weighted factorization of $x'$ in $M_i$. The left-weighted factorization of $\pphi_i(x')$ is given by $\Delta^k\pphi_i(s_1)\cdots \pphi_i(s_r)$. We have
\[\pphi_i(x')^s=\Delta^{k-1}\phi^{k-1}(\bbar{s})\pphi_i(s_1)\cdots\pphi_i(s_r)s.\]
If $s\notin\pphi_i(S_i)$, then the words $\phi^{k-1}(\bbar{s})\pphi_i(s_1)$ and $\pphi_i(s_r)s$ are in greedy normal form by Corollary \ref{cor:lcm_gcd_delta_product}. The expression above is then the left-weighted factorization of $\pphi_i(x')^s$ in $M$. Thus $\inf(\pphi_i(x')^s)=k-1$ and $\pphi_i(x')^s\notin\SSS(\pphi_i(x'))$. The same reasoning shows that if $\pphi_i(x')^{s^{-1}}\in \SSS(\pphi_i(x'))$, then $s\in \pphi_i(S_i)$. An immediate induction then shows that the connected component of $\pphi_i(x')$ in $\CG(\pphi_i(x))$ is contained in $\pphi_i(\SSS(x))$. Since $\CG(\pphi_i(x))$ is connected, we have $\SSS(\pphi_i(x))\subset \pphi_i(\SSS(x))$.

This also shows that the conjugacy graph $\CG(\pphi_i(x))$ is the image under of $\CG(x)$ under $\pphi_i$, whence the result on centralizers.
\end{proof}

\begin{prop}(Periodic elements in a $\Delta$-product)\label{prop:periodic_delta_product}\newline 
Let $p,q$ be nonzero integers, and let $\rho\in G(M)$ be a $(p,q)$-periodic element. There is some $i\in \intv{1,h}$ and some $(p,q)$-periodic element $\sigma\in G(M_i)$ such that $\rho$ is conjugate to $\pphi_i(\sigma)$ in $G(M)$.
\end{prop}
\begin{proof}
Let $\rho\in G(M)$ be a $(p,q)$-periodic element. If $\rho$ is conjugate to a power of $\Delta$, then the result is obvious. Otherwise, Proposition \ref{prop:periodic_elements} gives that $\rho$ is conjugate to an element of the form $\Delta^k s$ for some $s\in S$. By Proposition \ref{prop:simples_delta_product}, there is some $i\in \intv{1,h}$ and some $\ttilde{s}\in S$ such that $\pphi_i(\ttilde{s})=s$. We then have $\Delta^ks=\pphi_i(\Delta_i^k \ttilde{s})$. The element $\Delta_i^k \ttilde{s}$ is $(p,q)$-periodic in $M_i$ by Proposition \ref{prop:left-weighted_facto_in_delta_product}.
\end{proof}

Let $i\in \intv{1,h}$ and let $k_i$ be the smallest integer such that $\Delta_i^{k_i}$ is central in $G(M_i)$. By Proposition \ref{prop:simples_delta_product}, the smallest central power of $\Delta$ in $G(M)$ is given by the lcm of the $k_i$ for $i\in \intv{1,h}$.

\begin{prop}(Center of a nontrivial $\Delta$-product)\label{prop:centre_delta_product}\newline Assume that $h\geqslant 2$. The intersection in $G(M)$ of all the $\pphi_i(G(M_i))$ for $i\in \intv{1,h}$ is the subgroup generated by $\Delta$. The center of $G(M)$ is cyclic and generated by the smallest central power of $\Delta$ in $G(M)$.
\end{prop}
\begin{proof}
Since $h\geqslant 2$, we can choose $i\neq j$ in $\intv{1,h}$. Let $x_i\in G(M_i)$ and $x_j\in G(M_j)$ be such that $x:=\pphi_i(x_i)=\pphi_j(x_j)\in G(M)$. If $x_i=\Delta_i^{k_i}s_1\ldots s_r$ (resp. $x_j=\Delta_j^{k_j}t_1\ldots t_u$) is the left-weighted factorization of $x_i$ in $G(M_i)$ (resp. of $x_j$ in $G(M_j)$), then by Proposition \ref{prop:left-weighted_facto_in_delta_product}, the left-weighted factorization of $x$ in $G(M)$ is then given by 
\[x=\Delta^{k_i}\pphi_i(s_1)\cdots \pphi_i(s_r)=\Delta^{k_j}\pphi_j(t_1)\cdots \pphi_j(t_u).\]
We deduce that $k_i=k_j$,$r=u$ and $\pphi_i(s_k)=\pphi_j(t_k)$ for all $k\in \intv{1,r}$. Since $i\neq j$, Lemma \ref{lem:injectivity_of_pphi_i} gives that $r=0$ and $x$ is a power of $\Delta$.

Let now $x$ be an element of $Z(G(M))$, and let $i\in \intv{1,h}$. By definition, $x$ lies in the centralizer of $\pphi_i(a)$ for all $a\in A_i$. Since all the $M_i$ are distinct from $M(1,1)$, the elements of $A_i$ are not conjugate to a power of $\Delta$. Thus Proposition \ref{prop:sss_in_delta_product} gives that $x$ actually lies in the image under $\pphi_i$ of the common centralizer of all elements of $A_i$ in $M_i$. In other words we have $x\in \pphi_i(Z(G(M_i))$. We then have that $x\in \bigcap_{i=1}^h \pphi_i(Z(G(M_i)))$ is a power of $\Delta$, whence the result.
\end{proof}

\subsection{Hosohedral-type Garside groups}

\begin{definition} A Garside monoid $(M,\Delta)$ is a \nit{hosohedral-type monoid} if it is a $\Delta$-product of circular monoids. The enveloping group of a hosohedral-type monoid is a \nit{hosohedral-type group}.
\end{definition}

For instance, torus knot groups are hosohedral-type groups. Indeed, for $p,q$ two positive coprime integers, we have
\[\langle a,b~|~a^p=b^q\rangle\simeq G(M(1,p)*_\Delta M(1,q))\]

In his PhD thesis, Picantin introduced hosohedral-type monoids as the Garside monoids whose lattice of simple elements satisfy a strong property regarding gcds and lcms.

\begin{theo}\cite[Proposition 2.5]{thespicantin}\newline A finite lattice $(S,\wedge,\vee,0,1)$ has \nit{hosohedral-type} if any couple $(s,t)$ in $S$ satisfies
\[(a\wedge b,a\vee b)\in \{(a,b),(b,a),(0,1)\}\]
A Garside monoid $(M,\Delta)$ is a hosohedral-type monoid if and only if the lattice $(S,\preceq)$ of its simple elements has hosohedral-type. 
\end{theo}
Note that this theorem also covers the case of non-homogeneous Garside monoids, which we do not consider here.

In his proof, Picantin directly classifies all the Garside monoids whose lattice of simples is a fixed hosohedral-type lattice, in terms of the length of maximal chains. Under the notation of \cite[Proposition VI.2.5 and Definition VI.2.7]{thespicantin}, the monoid 
\[\mathrm{fus}^+\intv{h_1^{u_{1,1}}\cdot \ldots \cdot h_1^{u_{1,k_1}}\cdot h_2^{u_{2,1}}\cdot \ldots\cdot h_n^{u_{n,k_n}}}\]
is isomorphic to the $\Delta$-product
\[M(u_{1,1},h_1)*_\Delta\cdots *_\Delta M(u_{1,k_1},h_1)*_\Delta M(u_{2,1},h_2)*_\Delta\cdots *_\Delta M(u_{n,k_n},h_n).\]
This property regarding lcms and gcds of simple elements also induces strong geometric properties for the associated Garside group, as pointed out by Mireille Soergel in \cite[Theorem 4.6]{soergelsystolic}. By \cite[Proposition VI.1.11]{ddgkm}, every Garside monoid $(M,\Delta)$ with set of simples $S$ admits a presentation 
\begin{equation}\label{eq:garside_resentation} M\simeq \langle S~|~ s\cdot t=st~\forall s,t\in S\text{ such that }st\in S\rangle^+\end{equation}

\begin{theo}\cite[Theorem 4.6]{soergelsystolic}\newline Let $(M,\Delta)$ be a Garside monoid. The flag complex of the Cayley graph of $G(M)$ associated to presentation (\ref{eq:garside_resentation}) is systolic (in the sense of \cite[Section 2]{soergelsystolic}) if and only if $(M,\Delta)$ is a hosohedral-type monoid.
\end{theo}
Again this theorem also covers the case of non-homogeneous Garside monoids.

\begin{exemple}
In \cite[Question after Remark 4.8]{soergelsystolic}, Mireille Soergel considers the Garside group $G_k$, defined for an integer $k\geqslant 2$ by the presentation $G_k:=\langle a,b~|~ aba=b^k\rangle$. She asks whether or not this group is isomorphic to a hosohedral-type group. We have the following isomorphisms of groups, given by Tietze transformations
\begin{align*}
\langle a,b~|~aba=b^k \rangle&=\langle a,b,x~|~x=b^ka^{-1},aba=b^k \rangle\\
&=\langle a,b,x~|~a=x^{-1}b^k,aba=b^k \rangle\\
&=\langle b,x~|~x^{-1}b^kbx^{-1}b^k=b^k \rangle\\
&=\langle b,x~|~ b^{k+1}=x^2\rangle.
\end{align*}
The last group is a hosohedral-type group. More generally, a conjecture of Picantin (\cite[Conjecture 1]{thespicantin}) states that, if $M$ is a Garside monoid with two atoms, then the enveloping group $G(M)$ is isomorphic to either a group of the form $\langle a,b~|~a^p=b^q\rangle$ for positive integers $p,q$, or to an Artin group of dihedral type. In either case this would mean that the enveloping group of a Garside monoid with two generators is always a hosohedral-type group.
\end{exemple}

From now on, let $M=M(m_1,\ell_1)*_\Delta\cdots M(m_h,\ell_h)$ be a hosohedral-type monoid. Using Remark \ref{rem:identity_delta_product}, we assume that $(m_j,\ell_j)\neq (1,1)$ for all $j\in \intv{1,h}$. We also assume that $h\geqslant 2$, otherwise we recover results from Section \ref{sec:circular}. By Proposition \ref{prop:left-weighted_facto_in_delta_product}, we can identify the factors $M(m_j,\ell_j)$ (resp. $G(m_i,\ell_i)$) with the associated subgroup of $M$ (resp. of $G(M)$).

To avoid confusion, the simple $s(i,p)$ of the factor $M(m_j,\ell_j)$ will be denoted by $s_j(i,p)$.

\subsubsection{Conjugacy}

Like in the case of circular groups, an element in a super-summit set of a hosohedral-type monoid is either rigid or periodic.

\begin{prop}\label{prop:rigid_or_periodic_hosohedral}
Let $x$ be an element of $G(M)$. If $x$ lies in its own super-summit set, then it is either rigid or periodic.
\end{prop}

\begin{proof}
The proof imitates the case of circular groups.
First, if $\inf(x)=\sup(x)$, then we have $x=\Delta^k$ for some $k\in \Z$. In this case, $x$ is obviously both rigid and $(1,k)$-periodic.

Suppose now that $\sup(x)=\inf(x)+1$. We have $x=\Delta^k s_j(i,p)$ for some $j\in \intv{1,h}$, $i\in \intv{0,m_j-1}$ and $0<p<\ell_j$. The element $x$ is periodic if and only if it is periodic as an element of the factor $G(m_j,\ell_j)$. By Lemma \ref{lem:elements_periodiques_circulaires}, this is equivalent to $p+k\ell_j\equiv 0 [m_j]$. If this is not the case, the word $s_j(i,p)\phi^{-k}(s_j(i,p))$ is greedy in $M(m_i,\ell_i)$. It is then also greedy in $M$ by Proposition \ref{prop:left-weighted_facto_in_delta_product}.

Lastly, suppose that $\sup(x)>\inf(x)+1$. The left-weighted factorization of $x$ is given by $\Delta^ks_1\cdots s_r$ with $r>1$. We claim that $x$ is rigid. Otherwise, the word $s_r\phi^{-k}(s_1)$ is not greedy. By Corollary \ref{cor:lcm_gcd_delta_product}, this implies that $s_r$ and $s_1$ lie in the same factor $M(m_i,\ell_i)$. We can then apply the last part of the proof of Proposition \ref{prop:rigid_or_periodic} to show that we either have $\sup(\cyc(x))<\sup(x)$ or $\inf(\cyc(x))>\inf(x)$. In both cases, we have $x\notin\SSS(x)$.
\end{proof}

By Proposition \ref{prop:periodic_delta_product}, there are $(p,q)$-periodic elements in $G(M)$ if and only if there are $(p,q)$-periodic elements in some factor $G(m_j,\ell_j)$. However, two $(p,q)$-periodic elements coming from two different factors are not conjugate in general.

\begin{prop}\label{prop:periodic_hosohedral}
Let $p,q$ be integers. Two $(p,q)$-periodic elements of $G(M)$ are conjugate if and only if they both admit a conjugate lying the the same factor $G(m_j,\ell_j)$.
\end{prop}
\begin{proof}
Let $\rho,\sigma$ be two $(p,q)$-periodic elements in $G(M)$. If $\rho$ and $\sigma$ are conjugate, then $\SSS(\rho)=\SSS(\sigma)$. Let $x\in \SSS(\rho)$. We have either $x=\Delta^ks$ for some simple $s$ and some integer $k$ or $x=\Delta^k$. In both cases, $x$ is a conjugate of both $\rho$ and $\sigma$ lying in a fixed factor $M(m_i,\ell_i)$ by Proposition \ref{prop:simples_delta_product}.

Conversely, if $x$ and $y$ are respective conjugate of $\rho$ and $\sigma$ in $M(m_i,\ell_i)$, then $x$ and $y$ are $(p,q)$-periodic elements of $M(m_i,\ell_i)$. They are conjugate in $G(m_i,\ell_i)$ (in particular in $G(M)$) by Proposition \ref{prop:periodic_are_conjugate_circular}.
\end{proof}

\begin{exemple}\label{ex:periodic_torus}
Let $n$ be a positive integer, and consider $M=M(1,n)*_\Delta M(1,n)=\langle a,b~|~a^n=b^n\rangle^+$. By Proposition \ref{prop:sss_in_delta_product}, we have $\SSS(a)=\{a\}$ and $\SSS(b)=\{b\}$, thus $a$ and $b$ are two $(n,1)$-periodic elements that are not conjugate in $G(M)$.
\end{exemple}

On the other hand, the conjugacy of non-periodic elements in a hosohedral-type group behaves in the same way as in a circular group.

\begin{prop}\label{prop:sss_rigid_hosohedral}
Let $x\in G(M)$ be a non-periodic element. The super-summit set of $x$ is made of rigid elements. The only arrows starting from an object $y$ of $\CG(x)$ are labeled by $\init(y)$ and $\bbar{\fin(y)}$. In particular, one can go from an element of $\SSS(x)$ to any other by a finite sequence of cycling, decycling and application of the Garside automorphism.
\end{prop}
\begin{proof}
We mimic the proof of Proposition \ref{prop:conjugacy_graph_rigid_circular}. Let $y\in \SSS(x)$. Since $x$ is not periodic, $y$ is not periodic. It is then rigid by Proposition \ref{prop:rigid_or_periodic_hosohedral}, and we have $\sup(y)>\inf(y)$. We assume that the left-weighted factorization of $y$ is $\Delta^k s_{j_1}(i_1,p_1)\cdots s_{j_r}(i_r,p_r)$. Let $s$ be a simple element of $M$. We have
\[y^s=\Delta^{k-1}\phi^{k-1}(\bbar{s})s_{j_1}(i_1,p_1)\cdots s_{j_r}(i_r,p_r)s.\]
If $s\in \{1,\Delta\}$, then $y^s\in \SSS(x)$ is obvious. Otherwise, let $j$ be such that $s$ lies in the factor $M(m_j,\ell_j)$. If $j\notin \{j_1,j_r\}$, we have $\inf(y^s)=k-1$ and $y\notin \SSS(x)$.
\begin{enumerate}[-]
\item If $j_r=j_1$, then the proof of Proposition \ref{prop:conjugacy_graph_rigid_circular} gives that $y^s\in \SSS(x)$ if and only if $s\in \{\init(y),\bbar{\fin(y)}\}$.
\item If $j_r\neq j_1$ and $j=j_1$, then the word $s_{j_r}(i_r,p_r)s$ is greedy in $M$ by Corollary \ref{cor:lcm_gcd_delta_product}. We then have $y^s\in \SSS(x)$ if and only if $\phi^{k-1}(\bbar{s})s_{j_1}(i_1,p_r)=\Delta$, i.e. if $s=\init(y)$. 
\item If $j_r\neq j_1$ and $j=j_r$, then the word $\phi^{k-1}(\bbar{s})s_{j_1}(i_1,p_1)$ is greedy in $M$ by Corollary \ref{cor:lcm_gcd_delta_product}. We then have $y^s\in \SSS(x)$ if and only if $s_{j_r}(i_r,p_r)s=\Delta$, i.e. if $s=\bbar{\fin(y)}$.
\end{enumerate}
\end{proof}

\begin{theo}\label{theo:roots_conjugate_hosohedral}
Let $M$ be a hosohedral-type monoid. If $\alpha,\beta\in G(M)$ are such that $\alpha^n=\beta^n$ for some nonzero integer $n$, then we either have
\begin{enumerate}[-]
\item $\alpha$ and $\beta$ are conjugate.
\item $\alpha$ and $\beta$ are nonconjugate periodic elements of $G(M)$.
\end{enumerate}
\end{theo}
\begin{proof}
Again, $\alpha$ is periodic if and only if its power $\alpha^n$ is periodic, if and only if $\beta$ is periodic. Suppose that $\alpha$ and $\beta$ are not nonconjugate periodic elements of $G(M)$. Then $\alpha$ and $\beta$ are either conjugate periodic elements (in which case they are conjugate) or non-periodic elements. Up to replacing $\alpha,\beta$ with $\alpha^{-1},\beta^{-1}$, we now assume that $n>0$.

We assume from now on that $a$ and $b$ are non-periodic elements of $G(M)$. Up to conjugacy, we can assume that $a\in \SSS(a)$. By Proposition \ref{prop:rigid_or_periodic_hosohedral}, we have that $a$ is rigid. The element $x$ is then rigid as a power of the rigid element $a$. Let now $c\in G(M)$ be so that $b^c\in \SSS(b)$. Since $b$ is not periodic, $b^c$ is rigid as well as $x^c=(b^c)^n$. We have $x,x^c\in \SSS(x)$. By Proposition \ref{prop:sss_rigid_hosohedral}, there is an finite sequence of cycling, decycling, and applications of the Garside automorphism sending $x$ to $x^c$. By Lemma \ref{lem:cycling_power_rigid}, applying the same transformations to $a$ gives a rigid element $a'$ whose $n$-th power is $x^c$. Again by Lemma \ref{lem:cycling_power_rigid}, we have $a'=b^c$ and thus $a$ and $b$ are conjugate.
\end{proof}

\subsubsection{Partial computations of homology.}
Let $M=M(m_1,\ell_1)*_\Delta\cdots *_\Delta M(m_h,\ell_h)$ be a hosohedral-type monoid. We can construct the Dehornoy-Lafont complex (cf. Section \ref{sec:homology_circular}) associated to $M$ to try and compute the homology of $G(M)$. Since the lcm of two distinct atoms of a hosohedral-type group is $\Delta$ by Corollary \ref{cor:lcm_gcd_delta_product} and Lemma \ref{lem:gcd_lcm_simple_circulaire}, we can mimic the proof of Lemma \ref{lem:cellules_dehlaf_circular} to get

\begin{lem}\label{lem:cellules_dehlaf_hosohedral}
If $M$ is a hosohedral-type monoid with $m$ atoms, we order its atoms by $a_0<a_1<\ldots<a_{m-1}$. We have $\Xx_0=\{[\varnothing]\}$, $\Xx_1=\{[a_0],\ldots,[a_{m-1}]\}$, $\Xx_2=\{[a_0,a_i]~|~i\in \intv{1,m-1}\}$ and $\Xx_n=\varnothing$ for $n\geqslant 3$.
\end{lem}

Like in the case of a circular monoid, this implies that a hosohedral-type group has homological dimension at most $2$. Furthermore, as $H_2(G(M),\Z)$ is free, the computation of $G(M)^{\ab}$ is sufficient to determine all the integral homology of $G(M)$. Unfortunately, we are only able to compute the free part of $G(M)^{\ab}$.

\begin{prop}\label{prop:abelianisé_hosohedral}
Let $M=M(m_1,\ell_1)*_\Delta\cdots *_\Delta M(m_h,\ell_h)$ be a hosohedral-type monoid. For $i\in \intv{1,h}$, denote by $d_i$ the gcd of $m_i\wedge \ell_i$. The free part of $G(M)^{\ab}$ has rank $1+\sum_{i=1}^h (d_i-1)$.
\end{prop}
\begin{proof}
First, let $j\in \intv{1,h}$. Let also $e_0,\ldots,e_{d-1}$ denote the canonical basis of $\Z^{d_j}$. By Lemma \ref{lem:homology_circular}, the morphism $s_j(i,1)\mapsto e_{i'}$, where $i'$ is the remainder in the Euclidean division of $i$ by $d_j$ induces an isomorphism between $G(m_j,\ell_j)$ and $G(m,\ell)^{\ab}\simeq \Z^{d_j}$. Under this isomorphism, we have 
\[\Delta_j=s_j(0,m)\mapsto \left(\sum_{i=0}^{d_j-1}e_i\right)^{\frac{m_j}{d_j}}=\sum_{i=0}^{d_j-1}\frac{m_j}{d_j}e_i.\]
Let $A\simeq \Z^{\sum_{i=1}^n d_i}$ be the direct product of the $G(m_i,\ell_i)^{\ab}$ for $i\in \intv{1,h}$. The isomorphisms $G(m_j,\ell_j)^{\ab}\simeq \Z^{d_j}$ described above induce a morphism $p$ from the free product of the $G(m_i,\ell_i)$ to $A$, and we have $A=(G(m_1,\ell_1)*\cdots *G(m_h,\ell_h))^{\ab}$. We then have
\[G(M)^{\ab}=\bigslant{A}{\left\langle p(s_1(0,m_1))-p(s_i(0,m_i)),~\forall i\in \intv{2,h} \right\rangle}\]
as the vectors $p(s_1(0,m_1))-p(s_i(0,m_i))$ are linearly independent, they span a submodule of $A$ of rank $h-1$. The free part of $G(M)^{\ab}$ then has rank
\[\sum_{i=1}^h d_i-h+1=1+\sum_{i=1}^h(d_i-1).\]
\end{proof}

\begin{exemple}
Let $p,q\geqslant 2$ be integers, and consider the group 
\[G:=\langle a,b~|~a^p=b^q\rangle\simeq G(M(1,p)*_\Delta M(1,q)).\]
The proof of Proposition \ref{prop:abelianisé_hosohedral} gives that $G^{\ab}$ is isomorphic to the quotient of $\Z^2$ by the submodule generated by the vector $p(s_1(0,p))-p(s_2(0,q))=(-p,q)$. We obtain $H_1(G,\Z)=\Z\oplus \Z/d\Z$, where $d=p\wedge q$.
\end{exemple}

\printbibliography
\end{document}

%% file: pream.tex
\newcommand{\Xx}{\mathcal{X}}

\newcommand{\C}{\mathbb{C}}

\newcommand{\N}{\mathbb{N}}

\newcommand{\Z}{\mathbb{Z}}

\DeclareMathOperator\ho{Hom}

\DeclareMathOperator\GL{GL}

\DeclareMathOperator\cyc{cyc}
\DeclareMathOperator\dec{dec}
\DeclareMathOperator\SSS{SSS}

\DeclareMathOperator\md{md}

\DeclareMathOperator\ab{ab}

\DeclareMathOperator\CG{CG}
\DeclareMathOperator\init{init}
\DeclareMathOperator\fin{fin}

\newcommand{\bigslant}[2]{{\raisebox{.2em}{$#1$}\left/\raisebox{-.2em}{$#2$}\right.}}
\newcommand{\xdownarrow}[1]{%
  {\left\downarrow\vbox to #1{}\right.\kern-\nulldelimiterspace}
}

\newcommand{\infl}{\ar@{{(}->}}
\newcommand{\defl}{\ar@{->>}}

\newcommand{\intv}[1]{[\![#1]\!]}%Intervalles entiers
\newcommand{\muls}[1]{ \left\{\kern-0.6em \left\{ #1\right\}\kern-0.6em \right\} }%Multisets

\newcommand{\pphi}{\varphi}

\newcommand{\spa}{\vspace*{1ex}}

\newcommand{\bbar}[1]{\overline{#1}}

\newcommand{\ttilde}[1]{\widetilde{#1}}
\newcommand{\nit}[1]{\textbf{\emph{#1}}}